\documentclass{article}
\usepackage{amsmath,amsthm,amsfonts,amscd,eucal}

\numberwithin{equation}{section}

\def\ca{{\mathcal A}}
\def\cb{{\mathcal B}}
\def\cc{{\mathcal C}}
\def\cd{{\mathcal D}}
\def\ce{{\mathcal E}}

\def\ch{{\mathcal H}}
\def\cai{{\mathcal I}}

\def\ck{{\mathcal K}}
\def\cl{{\mathcal L}}
\def\cam{{\mathcal M}}

\def\car{{\mathcal R}}

\def\ct{{\mathcal T}}
\def\cu{{\mathcal U}}

\def\bc{{\mathbb C}}
\def\bg{{\mathbb G}}

\def\bn{{\mathbb N}}
\def\br{{\mathbb R}}

\def\a{\alpha}
\def\b{\beta}
        \def\G{\Gamma}
\def\d{\delta}        \def\D{\Delta}
\def\eps{\varepsilon}

\def\th{\vartheta}

\def\l{\lambda}       \def\La{\Lambda}
\def\m{\mu}
\def\n{\nu}

\def\p{\pi}
\def\r{\rho}
\def\s{\sigma}
    
\def\t{\tau}

\def\f{\varphi}
\def\c{\chi}
\def\om{\omega}        \def\O{\Omega}

\def\imply{\Rightarrow}

\def\ov{\overline}
\def\itm#1{\item{$(#1)$}}
\def\aff{\widehat\in}
\def\nat{\natural}
\def\supp{{\text{supp}}}
\def\ar{\ca^{\car}}
\def\au{\ca^{\cu}}
\def\at{\ca^{\ct}}
\def\e#1{{\rm e}^{#1}}
\DeclareMathOperator{\Lim}{Lim}

\def\bg{C$^{\infty}$-bounded geometry }
\def\tordim{{\mathfrak{tordim}}}

\def\inc{\nearrow}
\def\dec{\searrow}
\def\ri{\ov\car}
\def\rio{\ri_{0}}
\def\mt{\ov\cam}
\def\tlim{\t\,\text{lim }}

\def\pic#1{C_{0}(#1)}

\def\li{\ell_{0}^{\infty}(X)}
\def\ro{\car_{0}(X,\m)}
\def\Subset{\subset\subset}
\def\gN{{\mathfrak N}}                              

\def\sa{{\text{sa}}}

\def\ker{{\text{Ker}}}

\def\unia{{\ca^{**}}}

\newtheorem{Thm}{Theorem}[section]
\newtheorem{Cor}[Thm]{Corollary}
\newtheorem{Prop}[Thm]{Proposition}
\newtheorem{Lemma}[Thm]{Lemma}
\newtheorem{Sublemma}[Thm]{Sublemma}
\theoremstyle{definition}
\newtheorem{Dfn}[Thm]{Definition}
\newtheorem{exmp}[Thm]{Example}

\theoremstyle{remark}
\newtheorem{rem}[Thm]{Remark} 
\newtheorem{ack}{Acknowledgement} 
\title{\huge  Noncommutative Riemann integration\\
 and Novikov-Shubin invariants\\
 for open manifolds}
\author{
Daniele Guido, Tommaso Isola\\
Dipartimento di Matematica,\\ Universit\`a di Roma ``Tor
Vergata'',\\ I--00133 Roma, Italy.}

\date{}
\begin{document}
\maketitle
\markboth{Noncommutative Riemann integration}
{Noncommutative Riemann integration}
\renewcommand{\sectionmark}[1]{}
\bigskip

 \begin{abstract}
 Given a C$^*$-algebra $\ca$ with a semicontinuous semifinite trace 
 $\tau$ acting on the Hilbert space $\ch$, we define the family $\ar$ 
 of bounded Riemann measurable elements w.r.t. $\tau$ as a suitable 
 closure, {\it \`a la} Dedekind, of $\ca$, in analogy with one of the 
 classical characterizations of Riemann measurable functions 
 \cite{Tala1}, and show that $\ar$ is a C$^*$-algebra, and $\t$ 
 extends to a semicontinuous semifinite trace on $\ar$.

 Then, unbounded Riemann measurable operators are defined as the 
 closed operators on $\ch$ which are affiliated to $\ca''$ and can be 
 approximated in measure by operators in $\ar$, in analogy with 
 unbounded Riemann integration.  Unbounded Riemann measurable operators 
 form a $\tau$-a.e. bimodule on $\ar$, denoted by $\ov{\ar}$, and such 
 bimodule contains the functional calculi of selfadjoint elements of 
 $\ar$ under unbounded Riemann measurable functions.  Besides, $\t$ 
 extends to a bimodule trace on $\ov\ar$.

 Type II$_1$ singular traces for C$^{*}$-algebras can be defined on 
 the bimodule of unbounded Riemann-measurable operators.  
 Noncommutative Riemann integration, and singular traces for 
 C$^*$-algebras, are then used to define Novikov-Shubin numbers for 
 amenable open manifolds, show their invariance under 
 quasi-isometries, and prove that they are (noncommutative) asymptotic 
 dimensions.
 \end{abstract}

 \newpage

 \setcounter{section}{0}

 \section{Introduction.}\label{sec:zeroth}

In this paper we extend the notion of Riemann integrability to the non-abelian
setting, namely, given a C$^*$-algebra $\ca$ of operators acting on a Hilbert
space $\ch$, together with a semicontinuous semifinite trace $\t$, we define
both bounded and unbounded Riemann measurable operators on $\ch$ with respect to $\t$.
Then we apply the preceding construction in order to define the Novikov-Shubin
numbers for amenable open manifolds and show that they can be interpreted as asymptotic
dimensions.

\medskip

Bounded Riemann measurable operators associated with $(\ca,\t)$ form a 
C$^*$-algebra $\ar$, and $\t$ extends to a semicontinuous semifinite 
trace on it.  The unbounded Riemann measurable operators form a $\t$ 
a.e. $^*$-bimodule $\ov{\ar}$ on $\ar$ and $\t$ extends to a positive 
bimodule-trace on $\ov{\ar}$.  On the one hand, this leads to the 
notion of Riemann algebra as a C$^*$-algebra with some monotonic 
completion properties, and to the construction of the enveloping 
Riemann algebra for a pair $(\ca,\t)$, in analogy with what is done 
for the von~Neumann or Borel enveloping algebras \cite{Ped}.  On the 
other hand, Riemann integration provides a spatial theory of 
integration for C$^*$-algebras, namely the trace $\t$ extends to 
natural classes of bounded and closed unbounded operators on the given 
Hilbert space $\ch$.

The $\tau$-a.e. bimodule $\ov{\ar}$ of unbounded Riemann measurable 
operators is a natural environment for the definition of singular 
traces, in particular for the type II$_1$ singular traces considered 
in \cite{GI1}.  Indeed, while type I singular traces (or singular 
traces at $\infty$) may be defined on bounded operators (and are 
anyhow determined by their restriction to bounded operators, cf.  
\cite{GI1}), type II$_1$ singular traces (or singular traces at $0$) 
need unbounded operators to be defined, since they vanish on the 
bounded ones.  Therefore, we can define type 
II$_1$ singular traces for C$^*$-algebras by replacing the bimodule of unbounded 
measurable operators on a von Neumann algebra with the bimodule of 
unbounded Riemann measurable operators on a Riemann algebra.

As an application, and in fact a motivation, for our results on 
noncommutative Riemann integration, we study Novikov-Shubin numbers 
for amenable open manifolds.  In fact, the trace on almost-local 
operators on the manifold defined in \cite{GI2} can now be extended to 
Riemann measurable spectral projections.  Therefore the spectral 
density function of the $p$-Laplacian is well defined, and we may 
define the Novikov-Shubin numbers and show they are invariant under 
quasi-isometries.  As we shall see below, however, the noncommutative 
Riemann integration plays its major role in the dimensional 
interpretation of the Novikov-Shubin numbers.  Indeed such 
interpretation is related to the existence of singular traces (at 
zero) depending on the spectral asymptotics 
(at zero) of the $p$-Laplacian.  Unbounded Riemann integrable 
operators and singular traces on them furnish the necessary 
environment to prove such results.

\medskip

The first part of this paper concerns noncommutative Riemann integration.
The classical abstract Riemann integration (cf.  \cite{Kaplan,Tala1,Schw1}) is
based on a topological space $X$ (usually Hausdorff locally compact) and a
Radon measure $\mu$.  Riemann measurable functions may be defined in two
equivalent ways, either as functions which are discontinuous in a
$\mu$-negligible set or as functions which are $L^1$ approximated, both from
above and from below, by continuous functions.  The latter characterization may
be rephrased as follows: Riemann measurable functions are separating elements
of a Riemann-cut ($\car$-cut), where an $\car$-cut is a pair $(A^-,A^+)$ of
bounded classes of continuous functions in $C_0(X)$ s.t.
\begin{itemize}
	 \item $f^-\leq f^+$ for any $f^\pm\in A^\pm$
	 \item $\forall \eps>0$ $\exists f^\pm_\eps\in A^\pm : 
	 \int(f^+-f^-)d\mu<\eps$.
\end{itemize}
Replacing $(C_0(X),\mu)$ with $(\ca,\t)$, where $\ca$ acts on $\ch$, we define
the family $\ar$ of bounded Riemann measurable elements on $\ch$ w.r.t. $\t$ as
the linear span of the set of separating elements in $\ca''$ for $\car$-cuts in
$\ca$.  $\ar$ will also be called the $\car$-closure of $\ca$. Clearly $\t$
uniquely extends to a positive functional on $\ar$.

Then we prove that $\ar$ is a C$^*$-algebra, the extension of $\t$ to 
$\ar$ is a semicontinuous semifinite trace on it, $\ar$ is closed 
under functional calculus w.r.t. Riemann measurable functions, and is 
$\car$-closed, namely it is stable under $\car$-closure.  We call 
Riemann algebra an $\car$-closed C$^*$-algebra, and enveloping Riemann 
algebra of $\ca$ the $\car$-closure of $\ca$ in the universal representation.

The main technical problem concerning bounded Riemann measurable elements arises
in the definition already. Indeed we asked the classes $(A^-,A^+)$ of the
Riemann cut to be bounded. While such requirement is unnecessary in the abelian
case, due to the simple order structure there, it was apparently not known if
this was the case in the non-abelian setting either. The answer to this
question is negative in general, as it is shown in Subsection~\ref{sub:2.1} by
an explicit counterexample, therefore the two corresponding Riemann closures
(with or without boundedness assumption) are different in general.  Moreover
the non-uniformly bounded Riemann closure ($\cu$-closure) is a C$^*$-algebra
only when it coincides with the $\car$-closure. For these reasons we choose the
$\car$-closure in the definition of Riemann algebra.

Unbounded Riemann measurable elements w.r.t $\t$ are defined as closed operators
affiliated to $\ca''$ which may be approximated in measure by elements of
$\ar$.  This definition is inspired by the corresponding construction by
Christensen for Borel algebras \cite{Christensen}.  The family $\ov{\ar}$ of
such elements is a $\t$-a.e. $^*$-bimodule on $\ar$, namely it is closed under
the $^*$-bimodule operations and the $^*$-bimodule properties hold $\t$ almost
everywhere. Even if $\ov{\ar}$ is not necessarily closed w.r.t. the product
operation (see Remark \ref{noprod}), given two operators $S,\,T\in\ov{\ar}$,
there is $\widetilde T\in\ov{\ar}$, equal to $T$ $\t$-a.e., s.t. $S\widetilde
T\in\ov{\ar}$.  Moreover, $\t$ extends to a positive bimodule trace on it,
satisfying $\t(T^{*}eTe)=\t(eTeT^{*})$, for a suitable conull projection $e$.
Functional calculi of selfadjoint elements of $\ar$ under unbounded Riemann
measurable functions belong to $\ov{\ar}$.

\medskip

In the last part of this paper we extend the definition of the Novikov-Shubin
numbers to the case of amenable open manifolds.

As it is known, a general understanding of the geometric meaning of the
Novikov-Shubin invariants is still lacking. We believe that the definition of
these numbers as global invariants of an open manifold, rather then as homotopy
invariants of a compact one, may shed some light on their meaning. Our
interpretation of Novikov-Shubin numbers as asymptotic dimensions goes in this
direction.

This interpretation is based on a fundamental observation of Alain 
Connes, who showed that the integration on a $d$-dimensional compact 
Riemannian manifold may be reconstructed via the formula $\int f=\tau( 
f|D|^{-d} )$, where $\t$ is a singular trace, and $D$ is the Dirac 
operator.  Therefore, the dimension of a spectral triple in 
noncommutative geometry corresponds to the power of $|D|^{-1}$ giving 
rise to a non-trivial singular trace.  This is analogous to the  
situation of geometric measure theory, where the dimension determines 
which power of the radius of a ball gives rise to a non trivial volume 
on the space.

It was shown in \cite{GI5} that, as in the classical case, the Weyl 
asymptotics furnishes a dimension for a spectral triple via the 
formula $$d=\left(\limsup_{n\to\infty}\frac{\log\m_n}{\log 
1/n}\right)^{-1},$$ where $\m_n$ is the sequence of eigenvalues of 
$|D|^{-1}$, namely, when $d$ is finite non-zero, there exists a 
singular trace (not necessarily a logarithmic one) which is finite 
non-zero on $|D|^{-d}$.

This formula makes sense also for non-compact manifolds, if one 
replaces the eigenvalue sequence $\m_n$ with the eigenvalue function 
$\m(t)$ of $\D_{p}^{-1/2}$ (see Section 4), and recovers the dimension 
of the manifold.  But in this case, the behaviour for $t\to0$ may be 
considered too, giving rise to an asymptotic counterpart of the 
dimension.

Here we show that, under suitable assumptions, the asymptotic 
dimension associated with the Laplacian $\D_p$ on $p$-forms coincides 
with the $p$-th Novikov-Shubin number $\a_p$, and that it behaves as a 
dimension in noncommutative geometry, namely $\D_p^{-\a_p/2}$ gives 
rise to a non-trivial singular trace on the unbounded Riemann 
measurable operators, which is finite nonzero on $\D_p^{-\a_p/2}$.

The singular traceability of $\D_p^{-\a_p/2}$ has already been proved 
in \cite{GI5} in the case of $\G$-coverings.  The main problen in 
extending such dimensional interpretation to the case of amenable open 
manifolds is the fact that a normal trace on a von~Neumann algebra of 
$\G$-invariant operators is replaced by a semicontinuous trace on the 
C$^*$-algebra $\ca_{p}$ of almost local operators on $p$-forms defined in \cite{GI2}.  
Indeed, since Novikov-Shubin numbers are defined in terms of the 
spectral asymptotics near zero of the $p$-Laplacian, the needed 
singular traces should be looked for among type II$_1$ singular traces 
\cite{GI1}.  Such traces however are defined on bimodules of unbounded 
operators, since they vanish on bounded ones, and the notion of 
unbounded operator affiliated to a C$^*$-algebra is too restrictive 
for our purposes.  Unbounded Riemann measurable operators affiliated 
with $\ca_{p}$ and type II$_1$ singular traces on them furnish the 
environment for the traceability statements, hence for the dimensional 
interpretation.

However, since the semicontinuous trace on $\ca_{p}$ is not normal with respect to the 
given representation of $\ca_{p}$ on the space of $L^2$-differential 
forms, some assumptions are needed, as the vanishing of the torsion 
dimension introduced by Farber \cite{Farber}, cf. Remark \ref{r:remarks} 
$(c)$.

\section{Bounded Riemann integration.}\label{sec:first}

 Let $(\ca,\t)$ be a C$^*$-algebra with a semicontinuous semifinite 
 trace, acting on a Hilbert space $\ch$.  A pair $(A^-,A^+)$ of sets 
 in $\ca_\sa$ is called an $\car$-cut (w.r.t. $\t$) if $A^-,A^+$ are 
 uniformly bounded, {\it separated}, namely for any pair $a^{\pm}\in 
 A^\pm$ we have $a^-\leq a^+$, and $\t$-{\it contiguous}, namely
 \begin{equation}\label{eqn:contig}
	 \forall\, \eps>0\ \exists\, a_\eps^\pm\in A^\pm : \t(a^+_\eps-a^-_\eps)<\eps.
 \end{equation}
 An element $x\in\ca''$ is called separating for the $\car$-cut $(A^-, 
 A^+)$ if it is selfadjoint and for any $a^{\pm}\in A^\pm$ we have 
 $a^-\leq x\leq a^+$.

 \begin{Dfn} 
	 Let $(\ca,\t)$ be a C$^*$-algebra with a semicontinuous 
	 semifinite trace, and $\pi$ a faithful representation of $\ca$.  
	 Let us denote by $\car_{\pi}(\ca,\t)$ the linear span of the 
	 separating elements for the $\car$-cuts (w.r.t. $\t$) in 
	 $\pi(\ca)$.  When $\pi$ is the universal representation, we 
	 denote it simply by $\car(\ca,\t)$, and call it the {\it 
	 enveloping Riemann algebra} of $\ca$.  For the sake of 
	 convenience we use the shorthand notation $\ar\equiv 
	 \car_{\pi}(\ca,\t)$, when $\pi$ and $\t$ are clear from the 
	 context.  The C$^*$-algebra $\ca$ is called a {\it Riemann 
	 algebra} if it contains all the separating elements of its 
	 $\car$-cuts.
 \end{Dfn}

 If $x\in\ar_{\sa}$ and $(A^-,A^+)$ is a corresponding 
 $\car$-cut, either all or none of the $a$'s in $A^{\pm}$ have finite 
 trace.  In the first case we set $\t(x):=\inf \t(a^{+}) \equiv \sup 
 \t(a^{-})$.  In the second case, and if $x\geq0$, we set 
 $\t(x):=+\infty$.  This is the unique positive functional 
 extending $\t$ to $\ar$.

 The first property of a Riemann algebra $\ca$ is that it is closed 
 under Riemann functional calculus.  Moreover, the elements of $\ca$ 
 with $\t$-finite support are separating elements between upper and 
 lower Riemann sums made up with projections in $\ca$.

 \begin{Lemma}\label{Lemma:phieps} 
	 Let $(\ca,\t)$ be a C$^*$-algebra with a tracial weight, and 
	 denote by $\cai(\t)$ its domain.  Then, given a selfadjoint 
	 element $x$ in $\ca$,
	 $$
	 \t(\phi_\eps(|x|))<\infty, \forall\eps>0\Leftrightarrow x\in\ov{\cai(\t)},
	 $$
	 where, for any $\eps>0$, $\phi_\eps$ is an increasing continuous 
	 function from $\br_{+}$ to $\br$ such that $\phi_\eps=0$ on 
	 $[0,\eps/2)$ and $\phi_\eps=1$ on $(\eps,+\infty)$.
 \end{Lemma}

 \begin{proof} 
	 $(\Rightarrow)$ We have 
	 $|\t(x\phi_\eps(|x|))|\leq\|x\|\t(\phi_\eps(|x|))<\infty$, i.e. 
	 $x\phi_\eps(|x|)\in\cai(\t)$, and $\|x\phi_\eps(|x|)-x\|<\eps$, 
	 which implies $x\in\ov{\cai(\t)}$.  \\
	 $(\Leftarrow)$ By definition (\cite{Ped}, p.  175) 
	 $|x|\phi_\eps(|x|)$ belongs to the Pedersen ideal 
	 $K\left(\ov{\cai(\t)}\right)$.  Since $\cai(\t)$ is a dense ideal 
	 in $\ov{\cai(\t)}$, $K\left(\ov{\cai(\t)}\right)\subset\cai(\t)$ 
	 by minimality, hence $\t(|x|\phi_\eps(|x|))<\infty$.  Then 
	 $\t(\phi_\eps(|x|))<\frac2{\eps}\t(|x|\phi_\eps(|x|))<\infty$.
 \end{proof}

 Let us consider, for any selfadjoint $x$, the measure $\mu_x$ on 
 $\s(x)\setminus\{0\}$ defined by
 \begin{equation}\label{eqn:muex}
	 \int f(\l)d\m_x(\l)=\t(f(x)),\qquad f\in C_0(\s(x)\setminus\{0\}).
 \end{equation}

 \begin{Prop}\label{Prop:manyproj} 
	 Let $\ca$ be a Riemann algebra w.r.t. a semicontinuous semifinite 
	 trace $\t$.  Then: 
	 \item{$(i)$} 
	 $\ca$ is closed under Riemann functional calculus, namely
	 $$
	 \forall x \in \ca,\ f\in\car_0(\s(x)\setminus\{0\},\mu_x) \imply 
	 f(x)\in\ca, 
	 $$
	 where $\car_0(X,\mu)$ denotes the set of Riemann measurable 
	 functions on $X$ w.r.t. $\mu$ and vanishing at infinity (cf.  
	 Appendix).  In particular Riemann measurable spectral projections 
	 of selfadjoints elements of $\ca$ belong to $\ca$, where a 
	 spectral projection $e_{\O}(x)$ is Riemann measurable if 
	 $\O\Subset \s(x)\setminus\{0\}$ and $\m_{x}(\partial \O)=0$.  
	 \item{$(ii)$} 
	 Let $x\in\ca_{sa}$ have $\t$-finite support.  Then $x$ is the 
	 separating element between $\car$-cuts made by linear 
	 combinations of projections in $\ca$.  In particular, if $\t$ is 
	 densely defined on $\ca$, then $\ca$ is generated by its 
	 projections as a Riemann algebra.
 \end{Prop}

 \begin{proof} 
	 $(i)$.  It is sufficient to prove the assertion for a real-valued 
	 function $f$.  Then, by Proposition~\ref{Thm:1} in the Appendix, 
	 there exist an open set $V\subset\s(x)\setminus\{0\}$, with 
	 $\mu_x(V)$ finite, and functions $h$, $h^{\pm}_\eps\in 
	 C_0(\s(x)\setminus\{0\})$, with $h^{\pm}_\eps$ vanishing outside 
	 $V$, such that $\int (h^+_\eps-h^-_\eps) d\mu_x <\eps$ and 
	 $h^-_\eps\leq f-h\leq h^+_\eps$.  Then $h(x)$, 
	 $h^\pm_\eps(x)\in\ca$ by continuous functional calculus, and 
	 $f(x)-h(x)\in\ca$ because it is a Riemann algebra.  
	 \item{$(ii)$} 
	 Since $x$ has $\t$-finite support, the measure $\m_x$ is finite, 
	 therefore the set $S(x)=\{\a\in\s(x):\m_x(\a)\ne0\}$ is at most 
	 countable.  Consider the separated classes given by upper and 
	 lower Riemann sums of the function $f(t)=t$ on $\s(x)$, 
	 corresponding to subdivisions which do not intersect $S(x)$.  
	 Such classes are $\m_{x}$-contiguous and the corresponding 
	 spectral projections belong to $\ca$ by $(i)$.  The corresponding 
	 functional calculi of $x$ give the $\car$-cut for $x$.  
	 Concerning the last statement, we observe that $x\phi_\eps(|x|)$ 
	 has $\t$-finite support for any $\eps>0$ and any $x\in\ca_\sa$, 
	 by Lemma~\ref{Lemma:phieps}.  Then the thesis easily follows by 
	 part $(i)$.  
 \end{proof}

 Then we may state the main theorem of this section.

 \begin{Thm}\label{maintheorem} 
	 Let $\ca$ be a C$^*$-algebra with a semicontinuous semifinite 
	 trace $\t$, acting on a Hilbert space $\ch$.  Then 
	 \itm{i} $\ca^\car$ is a C$^*$-algebra 
	 \itm{ii} the above described extension of $\t$ to $\ar$ is a trace, and 
	 $(\ca^\car,\t)$ is a Riemann algebra 
	 \itm{iii} the GNS representation $\pi_\t$ of 
	 $\ca$ extends to a representation $\r_{\pi}$ of $\ca^\car$ into 
	 $\pi_\t(\ca)''$, and $\t|_{\ar}$ may be identified with the 
	 pull-back of the trace on $\p_{\t}(\ca)''$ via $\r_{\pi}$.  As a 
	 consequence $\t$ is semicontinuous semifinite on $\ar$.  
 \end{Thm}

 The rest of this section is devoted to the proof of 
 Theorem~\ref{maintheorem}.  In order to do that we have to introduce 
 {\it a priori} different kinds of Dedekind closures of $(\ca,\t)$, 
 namely we shall consider unbounded cuts ($\cu$-cuts), where the 
 uniform boundedness property is removed, and the corresponding 
 $\cu$-closure $\au$, and tight cuts ($\ct$-cuts), where the uniform 
 boundedness is strengthened by requiring that the $a^\pm_\eps$'s in 
 (\ref{eqn:contig}) verify
 $$
 \sup\s(a^+_\eps)<\sup\s(x)+\eps,\quad 
 \inf\s(a^-_\eps)>\inf\s(x)-\eps,
 $$
 and the corresponding $\ct$-closure, $\at$.  Of course we have 
 $\ca\subset\at\subset\ar\subset\au$, and $\t$ extends to $\au$ as 
 well.  We shall see that $\at=\ar$ is a C$^{*}$-algebra, while $\au$ 
 is not (cf. Example~\ref{notalg} below). It becomes a C$^{*}$-algebra 
 $iff$ it coincides with $\ar$.  
 This fact may be seen as a motivation for choosing $\ar$ instead of 
 $\au$ as the family of Riemann measurable elements.  As we shall see 
 in section 3, the possibility of taking products inside the set of 
 Riemann measurable elements is crucial for constructing $\ov{\ar}$.

 \begin{rem} In the abelian case, $\cu$-closure, $\car$-closure and
 $\ct$-closure coincide, indeed in this case we may always find a very tight
 cut, for which $\sup\s(a^+_\eps)=\sup\s(x)$ and $\inf\s(a^-_\eps)=\inf\s(x)$.  
 We conjecture that this is not always the case in the nonabelian setting.
 \end{rem}

 \bigskip

 \begin{Lemma}\label{Lemma:bimofirst} 
	 The sets of separating elements for $\car$-cuts and $\cu$-cuts 
	 are the selfadjoint parts of the $^*$-bimodules $\ar$ and $\au$ 
	 on $\ca$.
 \end{Lemma}

 \begin{proof} 
	 Linearity and $^*$-invariance are obvious.  Now let $a\in\ca$ and 
	 $x\in\ca''_{\sa}$ a separating element for the $\car$-cut 
	 $\{a^\pm_\eps\}$.  Then $a^*a^\pm_\eps a$ gives an $\car$-cut for 
	 $a^*xa$. \\
	 From the equalities
	 \begin{align} 
		 xy+yx	&=(x+1)^*y(x+1)-xyx-y\label{eqn:bimo}
		 \\
		 i(xy-yx)&=(x+i)^*y(x+i)-xyx-y\label{eqn:bimobis}
	 \end{align}
	 which hold for any pair of selfadjoint elements we then get the 
	 bimodule property for $\ca^\car$.  The proof for $\ca^\cu$ is 
	 analogous.
 \end{proof}

 \begin{Lemma}\label{def-rho}
 	There is a unique linear positive extension $\r:\au\to \pi_{\t}(\ca)''$ of the GNS 
 	representation of $\ca$.
 \end{Lemma}
 \begin{proof}
	Let $x\in\au_{\sa}$, and $\{a^{\pm}_{\eps}\}$ be a $\cu$-cut for 
	$x$.  Then $\pi_{\t}(a^{\pm}_{\eps})$ is a $\cu$-cut, and, 
	$\pi_{\t}(a^{+}_{\eps})$ and $\pi_{\t}(a^{-}_{\eps})$ converge to 
	the same element in $L^{1}(\pi_{\t}(\ca)'',\t)$, which is the 
	unique separating element of the $\cu$-cut, hence belongs to 
	$\pi_{\t}(\ca)''$.  Setting
	\begin{equation}\label{e:diagram}
		\r(x) := \lim_{\eps\to0}\pi_{\t}(a^{+}_{\eps}) = 
		\lim_{\eps\to0}\pi_{\t}(a^{-}_{\eps}),
	\end{equation}
	it follows easily that $\r(x)$ does not depend on the $\cu$-cut 
	and is linear and positive.
 \end{proof}
 
 \begin{Lemma} 
  \itm{i} $\r$ is a bimodule map. 
  \itm{ii} $\t|_{\au}=\t\circ\r$, hence it is a 
  trace on $\au$ as a bimodule on $\ca$.
 \end{Lemma}
 \begin{proof}
  $(i)$ Follows easily from equations (\ref{eqn:bimo}), 
  (\ref{eqn:bimobis}).  \\
  $(ii)$ Follows from positivity of $\r$ and $\t$.
 \end{proof}

 \begin{Lemma}\label{prop:auclosed}
  The $\cu$-closure of $\ca$ is $\cu$-closed, namely the set of 
  separating elements in $\ca_{\sa}''$ for $\cu$-cuts in $\au$ is 
  contained in $\au$.
 \end{Lemma}
 \begin{proof}
  Indeed if $\{x^\pm_{\eps'}\}_{\eps'>0}$ is a $\cu$-cut in $\au$ for 
  $x\in \ca_\sa''$, and $\{a(x^\pm_{\eps'})^\pm_\eps\}_{\eps>0}$ is a 
  $\cu$-cut in $\ca$ for $x^\pm_{\eps'}$, then 
  $\{(a(x^+_{\eps})^+_\eps,a(x^-_{\eps})^-_\eps)\}$ gives a $\cu$-cut 
  in $\ca$ for $x$.
 \end{proof}
 
 With the above terminology, Theorem \ref{maintheorem} shows that the 
 $\car$-closure of a C$^*$-algebra is $\car$-closed.
  
 \begin{Lemma}\label{Lemma:contfunctcalc}
	 Let $\ca$ be a C$^*$-algebra acting on a Hilbert space $\ch$, 
	 with a semicontinuous semifinite trace $\t$.  If 
	 $C_0(x):=\{f(x):f\in C_0(\s(x)\setminus\{0\})\}$ is contained in 
	 $\ca^\cu$ for a selfadjoint element $x$, then $C_0(x)\subset\at$.  
	 In particular $\at$, $\ar$ and $\au$ have the same projections, 
	 and any C$^{*}$-subalgebra of $\au$ is actually contained 
	 in $\at$.
 \end{Lemma}
 \begin{proof}
	First we observe that any projection $e$ in $\ca^\cu$ belongs to 
	$\ct$.  For any $\d>0$ we consider the operator increasing 
	functions (cf.  \cite{Ped})
    $$
    f^+_\d(z)=\frac{(1+\d)z}{\d+z}\quad;
    \qquad f^-_\d(z)=\frac{\d z}{1+\d-z}\quad.
    $$
	If $\{a^\pm_\eps\}_{\eps>0}$ gives a $\cu$-cut for $e$, 
	$\{f^\pm_\d(a^\pm_\eps)\}_{\eps>0}$ gives an $\car$-cut for any 
	$\d>0$, indeed $f^\pm_\d(e)=e$ together with operator monotonicity 
	imply that $f^-_\d(a^-_\eps)\leq e\leq f^+_\d(a^+_\eps)$, and we 
	have the estimate
    \begin{align*}
	   \t(f^+_\d(a^+_\eps)-&f^-_\d(a^-_\eps))\\
       =&\t(f^+_\d(a^+_\eps)-f^+_\d(e))+\t(f^-_\d(e)-f^-_\d(a^-_\eps))\\
       \leq&(1+\d)\d [\t( (a^+_\eps+\d)^{-1}(a^+_\eps-e)(e+\d)^{-1})+\\
       &\t( (1+\d-e)^{-1} (e-a^-_\eps)(1+\d-a^-_\eps)^{-1})]\\
       \leq&\frac{1+\d}{\d}\eps.
    \end{align*}
	Moreover, we have $\sup\s(f^+_\d(a^+_\eps))\leq1+\d$ and 
	$\inf\s(f^-_\d(a^-_\eps))\geq-\d$, therefore we may easily extract 
	a $\ct$-cut for $e$.

	Now let $C_0(x)\subset\cu$; we want to show that $x\in\at$.  First 
	observe that $C_0(x^\pm)\subset\cu$, where $x^\pm$ is the 
	positive, resp.  negative part of $x$.  If $0$ is in the convex 
	hull of the spectrum of $x$, then $\inf\s(x^\pm)=0$, hence 
	$\ct$-cuts for $x^\pm$ give a $\ct$-cut for $x$, namely we may 
	restrict to the case $\inf\s(x)=0$.  If $\inf\s(x)=m>0$, then 
	$I\in\au$ by hypothesis, and it is sufficent to prove the 
	statement for $x-mI$, while if $\sup\s(x)=m<0$, it is sufficent to 
	prove the statement for $-x+mI$, therefore again we may suppose 
	$\inf\s(x)=0$.  Hence, by multiplying with an appropriate costant, 
	we may assume $\inf\s(x)=0$ and $\sup\s(x)=1$.
 
	Then we fix a $\d>0$ and apply Lemma~\ref{Lemma:abelian2} in the 
	Appendix to the identity function on the spectrum of $x$.  We may 
	therefore find a sequence $e_n$ of Riemann measurable spectral 
	projections of $x$ and a sequence of positive numbers $\a_n$ such 
	that $\sum_n\a_n e_n=x$ and $\sum_n \a_n=1$.

	Moreover, since $\ca^\cu$ is $\cu$-closed by Proposition 
	\ref{prop:auclosed}, all the Riemann measurable spectral 
	projections of $x$ belong to $\ca^\cu$ (cf.  the proof of 
	Proposition~\ref{Prop:manyproj}). 

	Therefore, for any $e_n$ we obtain an $\car$-cut 
	$(\{a^-_{n,\eps}\}_{\eps>0},\{a^+_{n,\eps}\}_{\eps>0})$ such that 
	$\sup_{\eps}\s(a^+_{n,\eps})\leq1+\d$ and 
	$\inf_{\eps}\s(a^-_{n,\eps})\geq-\d$.  Hence, setting 
	$a^\pm_\eps=\sum_n\a_na^\pm_{n,\eps}$, we obtain
    \begin{equation*}
    	a^-_\eps\leq x\leq a^+_\eps,\cr 
		\t(a^+_\eps-a^-_\eps)\leq\sum_n\a_n
		\t(a^+_{n,\eps}-a^-_{n,\eps})\leq\eps, \cr
		\sup\s(a^+_{\eps})\leq\sum_n\a_n
		\sup\s(a^+_{n,\eps})\leq 1+\d, \cr
		\inf\s(a^-_\eps)\geq\sum_n\a_n\inf\s(a^-_{n,\eps})\geq-\d.
    \end{equation*}
	From these $\car$-cuts it is easy to extract a $\ct$-cut.  
	Finally, we observe that $C_0(f(x))\subset\au$ for any $f\in 
	C_0(\s(x)\setminus\{0\})$, hence $C_0(x)\subset\at$.
 \end{proof}

 \begin{Lemma}\label{Lemma:normclosure}
 	Let $X\subset\at$ be a $^*$-invariant vector space. Then $\ov{X}\subset\ar$.
 \end{Lemma}
 \begin{proof}
	Let $x\in \ov{X}_\sa$, $x_n\to x$ in norm, $x_n\in X_\sa$.  First, 
	possibly passing to a subsequence, we may assume 
	$\|x-x_n\|<2^{-n}$.  Then set $y_n=x_{n+1}-x_n\in X$ and observe 
	that $\|y_n\|<3\cdot 2^{-(n+1)}$.  By hypothesis we may find a  
	$\ct$-cut in $\ca$ for $y_n$, hence elements 
	$a^\pm_{n\eps}\in\ca$ such that $a^-_{n\eps}\leq y_n\leq 
	a^+_{n\eps}$, $\|a^\pm_{n,\eps}\|\leq2\cdot 2^{-n}$, 
	$\t(a^+_{n\eps}-a^-_{n\eps})\leq\eps\ 2^{-n}$.  Then 
	$a^\pm_\eps:=\sum_n a^\pm_{n\eps}$ belongs to $\ca$ by uniform 
	convergence, and 
	$$
    \t(a^+_\eps-a^-_\eps)\leq\sum_n\t(a^+_\eps-a^-_\eps)\leq\eps,
    $$
	where the first inequality follows by the semicontinuity of $\t$.  
	Since $\|a^\pm_\eps\|\leq\sum_n\|a^\pm_{n\eps}\|\leq2$, we get 
	$x\in\ca^\car$.
 \end{proof}
 
 \begin{Lemma}\label{Lemma:resolvent}
	 Let $x\in\ar$, and consider the functions $f_{t}(z):= 
	 \frac{tz}{1-tz}$, $t\in[0,1)$.  Then $f_{t}(x)\in\ar$, and 
	 $\r(f_{t}(x))=f_{t}(\r(x))$, for sufficiently small $t$.
 \end{Lemma}
 \begin{proof}
	 Let $x$ be a separating element for the $\car$-cut 
	 $\{a^\pm_\eps\}\subset\ca$, with $r=\sup_\eps\|a^\pm_\eps\|$.  
	 Then $f_{t}(a^\pm_\eps)\in\ca$ for any $t\in[0,1/r)$.  Since, for 
	 any such $t$, the function $f_{t}$ is operator monotone 
	 increasing (cf.  \cite{Ped}) on $(-\infty,r)$, we get,
     $$
     f_{t}(a^-_\eps)\leq f_{t}(x)
     \leq f_{t}(a^+_\eps) \qquad t\in[0,1/r)
     $$
     and
     \begin{align*}
		 f_{t}(a^+_\eps)-f_{t}(a^-_\eps)
		 &=(I-tx^+_\eps)^{-1}(tx^+_\eps(I-tx^-_\eps)-(I-tx^+_\eps)tx^-_\eps) 
		 (I-tx^-_\eps)^{-1}\cr 
		 &=t(I-tx^+_\eps)^{-1}(x^+_\eps-x^-_\eps)(I-tx^-_\eps)^{-1}.
	 \end{align*}
     Then, taking $0\leq t<\frac1{2r}$, we get 
     \begin{align*}
		 \left\| f_{t}(a^\pm_\eps) \right\| & \leq 1  \cr
		 \t\left( 
		 f_{t}(a^+_\eps)-f_{t}(a^-_\eps) 
		 \right) & \leq\frac{2\eps}{r},
	 \end{align*}
	 which means that $f_{t}(x)$ is a separating element for the 
	 $\car$-cut $\{f_{t}(a^\pm_\eps)\}$, therefore, $f_{t}(x)$ belongs 
	 to $\ar$ for any $t\in[0,\frac1{2r})$.  Finally, making use of 
	 equation~(\ref{e:diagram}) (where the limits are in $L^1$) we get
	 $$
	 \r(f_{t}(x)) = \lim_{\eps\to0} \pi_{\t}(f_{t}(a_{\eps}^{\pm})) 
	 = \lim_{\eps\to0} f_{t}(\pi_{\t}(a_{\eps}^{\pm})) = f_{t}(\r(x)).
	 $$
 \end{proof}

 \begin{Lemma}\label{keylemma} 
	 Let $\ca$ be a C$^*$-algebra acting on a Hilbert space $\ch$, 
	 with a semicontinuous semifinite trace $\t$.  If the norm closure 
	 of $\ar$ is contained in $\ca^\cu$, then $\at=\ar$, and 
	 Theorem~\ref{maintheorem} holds for $\ca$.
 \end{Lemma}
 \begin{proof} 
	 The proof requires some intermediate steps.  Some of the 
	 arguments are taken form the proof of Kadison \cite{Kadi1} that 
	 the Borel closure of a C$^*$-algebra is a C$^*$-algebra.  \\
 Step $(i)$. $\at=\ar$ is a Banach space and is closed under
 continuous functional calculus. 
	
	 Let $x\in\ar$, then by Lemma~\ref{Lemma:resolvent} we have 
	 $f_{t}(x)\in\ar$, and the equality
	 \begin{equation}\label{eq:powers}
		 x^n=\lim_{t\to0}\left( f_{t}(x)-\sum_{k=1}^{n-1} t^k 
		 x^k\right)t^{-n},
	 \end{equation}
	 inductively implies $x^n\in\au$, hence $C_0(x)\subset\au$.  Then 
	 Lemma~\ref{Lemma:contfunctcalc} shows that $C_0(x)\subset\at$, 
	 namely $\at=\ar$ and it is closed under functional calculus.  
	 Then, since $\ar$ is a vector space, we may apply 
	 Lemma~\ref{Lemma:normclosure} to $\at$, obtaining that $\ar$ is 
	 norm closed.  \\
 Step $(ii)$.	$\ar$ is $\car$-closed, namely if 
	 $x\in \au$ is a separating element for an $\car$-cut 
	 $\{x_{\eps}^{\pm}\}\subset\ar$, then $x\in\ar$.  

	 Choose a $\ct$-cut $\{(a_{\eps}^{\pm})_{\d}^{\pm}\}\subset \ca$ 
	 for $x_{\eps}^{\pm}$.  Then 
	 $\{(a_{\eps}^-)_{\eps}^-,(a_{\eps}^+)_{\eps}^+\}$ is an 
	 $\car$-cut for $x$, and $x\in\ar$.  \\
 Step $(iii)$.	$\r$ is norm-continuous and $\r(x^2)=\r(x)^2$. 
  
	 Indeed, adjoining the identity to $\ca$ if necessary, if 
	 $x\in\ar_{\sa}$, then, as $\r$ is a linear and positive map, we 
	 get $\|\r(x)\|\leq \|x\| \|\p(1)\|$, and by linearity of $\r$ we 
	 get norm continuity.  The last equation follows by applying $\r$ 
	 to equation~(\ref{eq:powers}) for $n=2$ and using 
	 Lemma~\ref{Lemma:resolvent} \\
 Step $(iv)$. If $x,y\in\ar_\sa$, then
	 $xy+yx\in\ar$ and $\r(xy+yx)=\r(x)\r(y)+\r(y)\r(x)$.

	 Immediate by steps $(i)$ and $(iii)$, and 
	 \begin{equation}\label{eqn:quadrato}
	 	xy+yx=(x+y)^2-x^2-y^2.
	 \end{equation}
	 \\
 Step $(v)$.	If $x,y\in\ar_\sa$, then $xyx\in\ar$ and $\r(xyx)=\r(x)\r(y)\r(x)$.

     Immediate by step $(iv)$ and 
	 \begin{equation}\label{eqn:conj}
		 2xyx=(xy+yx)x+x(xy+yx)-(yx^2+x^2y).
	 \end{equation}
	 \\
 Step $(vi)$.	 If $x\in\ar_\sa$, $y\in\ca_\sa$ then 
	 $(x+i)^*y(x+i)\in\ar$ and 
	 $\r((x+i)^*y(x+i))=(\r(x)+i)^*\r(y)(\r(x)+i)$.  

	 Follows by equations (\ref{eqn:bimobis}) and (\ref{eqn:conj}) for 
	 $x\in\ar_\sa$, $y\in\ca_\sa$, and by the bimodule property.  
	 \\
 Step $(vii)$. If $x,y\in\ar_\sa$,  then 
	 $(x+i)^*y(x+i)\in\ar$ and 
	 $\r((x+i)^*y(x+i))=(\r(x)+i)^*\r(y)(\r(x)+i)$.

	 Let $\{b^{\pm}_{\eps}\}$ be an $\car$-cut for $y$, then we have 
	 $$(x+i)^{*}b_{\eps}^{-}(x+i)\leq (x+i)^{*}y(x+i) \leq
	 (x+i)^{*}b_{\eps}^{+}(x+i),$$
	 and, by step $(vi)$,
	 \begin{align*}
		 \t((x+i)^{*}(b_{\eps}^{+}-b_{\eps}^{-})(x+i)) 
		 & = \t\circ\r((x+i)^{*}(b_{\eps}^{+}-b_{\eps}^{-})(x+i)) \\
		 & \leq \|x+i\|^{2}\t(b_{\eps}^{+}-b_{\eps}^{-}) \leq \eps \|x+i\|^{2}
	 \end{align*}
	 so that $\{(x+i)^{*}b_{\eps}^{\pm}(x+i)\}\subset\ar$ is an 
	 $\car$-cut for $(x+i)^{*}y(x+i)$.  Then the thesis follows by 
	 step $(ii)$ and equation~(\ref{e:diagram}).  \\
 Step $(viii)$. $\ar$ is a C$^*$-algebra and $\r$ is a representation.

	 The only missing property is multiplicativity, which directly 
	 follows from steps $(iv)$ and $(vii)$, and 
	 equation~(\ref{eqn:bimobis}).  \\
 Step $(ix)$. $\r$ is the GNS representation of $\ar$ w.r.t. $\t$.

	 Let us denote by $(\pi_{\t},\ch_{\t},\eta_{\t})$ the GNS triple of 
	 $\ca$ w.r.t. $\t$.  Set $\gN:= \{x\in\ar:\t(x^{*}x)<\infty\}$, 
	 and define the linear map $\eta:\gN\to\ch_{\t}$ as follows.  If 
	 $x\in\gN_{+}$, and $\{a_{\eps}\}\subset\gN\cap\ca_{+}$ is s.t. 
	 $x\leq a_{\eps}$ and $\t(|a_{\eps}-x|^{2})<\eps$, let us set 
	 $\eta(x):= \lim_{\eps\to0}\eta_{\t}(a_{\eps})$, where the limit 
	 is  independent of $\{a_{\eps}\}$, and extend 
	 by linearity to all of $\gN$, which is generated by its positive 
	 elements.  Then it is easy to see that 
	 $(\eta(x),\eta(y))=\t(x^{*}y)$, for $x,y\in\gN$, and all that 
	 remains to show is that $\r(x)\eta(y)=\eta(xy)$ for $x\in\ar$, 
	 $y\in\gN$, and it sufficies to show it for $x,y\geq0$.  Let us 
	 prove it first for $y\in\gN\cap\ca_{+}$.  Observe that, using 
	 equations (\ref{eqn:bimo}) and (\ref{eqn:bimobis}), we can write
	 \begin{align}\label{e:xy}
	 	xy & = \frac12 (y+1)x(y+1) + \frac{i}{2} (y+i)^{*}x(y+i) -\frac{1+i}{2} 
	 	yxy-\frac{1+i}{2}x \notag\\
		& \equiv \sum_{k=1}^{4}\l_{k}b_{k}^{*}xb_{k},
	 \end{align}
	 where $\l_{k}\in\bc$, $b_{k}\in\ca+\bc$. Let now 
	 $\{a^{\pm}_{\eps}\}\subset\ca$ be an $\car$-cut, with $x$ as 
	 separating element. Then, using equation (\ref{e:xy}), we get
	 	$a^{+}_{\eps}y= \sum_{k=1}^{4}\l_{k}b_{k}^{*}a^{+}_{\eps}b_{k}$,
	 so that $b_{k}^{*}xb_{k}\leq b_{k}^{*}a^{+}_{\eps}b_{k}$, and 
	 \begin{equation*}
	 	\t(|b_{k}^{*}(a^{+}_{\eps}-x)b_{k}|^{2})\leq 
		\|b_{k}\|^{2}\|a^{+}_{\eps}-x\| \t(a^{+}_{\eps}-x)\to 0,\ \eps\to0,
	 \end{equation*}
	 hence,
	 \begin{align*}
		\lim_{\eps\to0} \pi_{\t}(a^{+}_{\eps})\eta_{\t}(y) & =
		\lim_{\eps\to0} \eta_{\t}(a^{+}_{\eps}y) =
		\lim_{\eps\to0}\sum_{k=1}^{4}\l_{k}\eta(b_{k}^{*}a^{+}_{\eps}b_{k}) 
		\\
		& = \sum_{k=1}^{4}\l_{k}\eta(b_{k}^{*}xb_{k}) = \eta(xy). 
	 \end{align*}
	 But we also have
	 \begin{equation*}
		 \|[\pi_{\t}(a^{+}_{\eps})-\r(x)]\eta_{\t}(y)\|^{2} \leq
		 \|\pi_{\t}(a^{+}_{\eps})-\r(x)\|\, 
		 \|\pi_{\t}(a^{+}_{\eps})-\r(x)\|_{1}\, \|y\|^{2} \to0,
	 \end{equation*}
	 as $\eps\to0$, so that $\r(x)\eta(y)=\eta(xy)$. Finally for 
	 $x\in\ar_{+}$, $y\in\gN_{+}$, $a\in\gN\cap\ca_{+}$, we get
	 \begin{equation*}
	 	(\eta(a),\r(x)\eta(y)) = (\r(x)\eta(a),\eta(y)) = 
	 	(\eta(xa),\eta(y)) = \t(axy) = (\eta(a),\eta(xy))
	 \end{equation*}
	 and the thesis follows.
 \end{proof}

 In the following subsections we shall prove that the closure of
 $\ar$ is contained in $\au$ in the cases in which $\ca$ is unital and 
 the trace is finite, $\ca$ is non-unital and the trace is densely 
 defined, and finally in the non-densely defined case, thus 
 completing the proof of Theorem~\ref{maintheorem}.
 
 \begin{Cor} 
	 $\ar$ is the unique Riemann algebra between $\ca$ and $\au$.  
	 $\au$ is a C$^*$-algebra if and only if it coincides with $\ar$.
 \end{Cor}
 \begin{proof} 
	By Lemma \ref{keylemma} and Theorem \ref{maintheorem} we obtain 
	that $\ca^{\car}$, being the minimal Riemann algebra between $\ca$ 
	and $\ca^{\cu}$ and being the maximal C$^{*}$-algebra between 
	$\ca$ and $\ca^{\cu}$, is indeed the unique Riemann algebra between 
	$\ca$ and $\au$.  As a consequence, either $\ar=\au$, or $\au$ is 
	not a C$^{*}$-algebra.
 \end{proof}
 
 \begin{rem}
	We say that a linear subspace $X$ in $\unia$ is completely 
	$\t$-measurable if, for any faithful representation $\pi$ of 
	$\ca$, the map 
	$\r_\pi:=\pi_\t\cdot\pi^{-1}:\pi(\ca)\to\pi_{\t}(\ca)$ extends to 
	a linear map from $\pi(X)$ in such a way that the following 
	diagram is commutative.
    \begin{equation*}
    \begin{CD}
    X		@>>> 	 \unia		\\
    @V{\pi}VV       @VV{\pi_{\t}}V	\\
    \pi(X)		@>{\rho_{\pi}}>> \pi_\t(\ca)''	
    \end{CD}
    \end{equation*}
	When $\pi$ is not normal, $\t$ cannot be extended to $\pi(\ca)''$.  
	However, if $\car$ is a completely measurable C$^{*}$-algebra, 
	$\ca\subset \car\subset \unia$, $\t$ uniquely extends to a trace 
	on $\pi(\car)$.  One can show that $\cu(\ca,\t)$, hence 
	$\car(\ca,\t)$, are completely measurable, however the algebra 
	generated by $\cu(\ca,\t)$ may fail to be completely measurable.  
	Observe that $\pi(\cu(\ca,\t))\subset \au$ and 
	$\pi(\car(\ca,\t))\subset \ar$, whereas it is not known to the 
	authors if equalities hold.
 \end{rem}


\subsection{Finite trace on a unital C$^*$-algebra}\label{sub:2.1}

 In this subsection we prove Theorem~\ref{maintheorem} in the case of 
 a unital C$^*$-algebra with a finite trace.  This will be an 
 immediate corollary of the following theorem.

 \begin{Thm}\label{Thm:1main} Let $(\ca,\t)$ be a unital 
 C$^{*}$-algebra acting on a Hilbert space $\ch$, $\t$ a tracial 
 state on $\ca$, then $\au$ is norm closed, therefore
 Theorem~\ref{maintheorem} holds for $\ca$.
 \end{Thm}
 \begin{proof}
	Let $x\in\ov{\au}$ and $x_{n}\in\au$, $\|x_{n}-x\|<\frac{1}{n}$.  
	Then $x$ is a separating element for the $\cu$-cut 
	$(\{x_{n}-\frac{1}{n}\}, \{x_{n}+\frac{1}{n}\})$, so that 
	$x\in\au$ since $\au$ is $\cu$-closed.  Then the result follows by 
	Lemma~\ref{keylemma}.
 \end{proof}

 Since $\au$ is a norm closed $^{*}$-bimodule, $\ar\neq\au$ $\iff$ 
 $\au$ is not an algebra.  We conclude this subsection with an example 
 showing that this is indeed possible.  Apparently such phenomenon 
 depends on $\t$ being non-faithful on $\ca^\car$, which can happen 
 even if $\t$ is faithful on $\ca$.

 \begin{exmp}\label{notalg} Let $\mu$ be the sum of the Lebesgue measure on
 $[0,1]$ plus the Dirac measure in $\{0\}$, and consider the C$^*$-algebra 
 $\ca=\{f\in C([0,1])\otimes M(2):f_{12}(0)=f_{21}(0)=0\}$ acting on the 
 Hilbert space $\ch=L^2([0,1],d\mu)\otimes\bc^2$.  Let $\t$ be the 
 state on $\cb(\ch)$ defined by
 $$
 \t(f):=\frac12\int \mbox{tr}f(t) dt+\frac12f_{11}(0).
 $$
 Then $\t$ is a trace on $\ca''$ and is faithful on $\ca$, $\ca^\car$ 
 is given by the $M(2)$-valued Riemann measurable functions $f$ on 
 $[0,1]$ s.t.  $f_{11}$ is continuous in $0$ and $f_{12}$, $f_{21}$ 
 vanish and are continuous in $0$, $\ca^\cu$ is given by the 
 $M(2)$-valued Riemann measurable functions $f$ on $[0,1]$ s.t.  
 $f_{11}$ is continuous in $0$ and $f_{12}$, $f_{21}$ vanish in $0$.  
 In particular $\ca^\cu$ is not an algebra.
 \end{exmp}

 \begin{proof} 
	 Since $\ca''$ is given by $\{f\in L^\infty([0,1],d\m)\otimes M(2) 
	 : f_{12}(0)=f_{21}(0)=0\}$, the trace property follows.

 Since elements in $\ca^\cu$ are separating elements for 
 $\cu$-cuts in $\ca$, we obtain that $f_{11}$ is continuous in $0$.  
 Since $\ca^\car$ is a $^*$-algebra, $f\in\ca^\car$ 
 implies $|f|^2\in\ca^\car$, hence $|f_{11}|^2+|f_{12}|^2$ is 
 continuous in $0$.  Since $f_{12}$ vanishes in $0$, $f_{12}$ (and 
 analogously $f_{21}$) is continuous in $0$.  On the other hand, let 
 $g$ be a real valued Riemann measurable function on $[0,1]$ w.r.t. 
 Lebesgue measure s.t.  
 $g(0)=0$ and $|g|\leq1$ and consider the matrix 
 $f:=\left(\begin{matrix}0&g\cr g&0\cr\end{matrix}\right)$.  For any 
 $\eps>0$ set
 $$
 f^\pm_\eps(t):=
 \begin{cases}
 \left(
 \begin{matrix}
 0&g(t)\cr 
 g(t)&0\cr
 \end{matrix}
 \right)                & t>\eps \cr
 \cr
 \left(
 \begin{matrix}
 \pm\eps^2&0\cr 
 0&\pm\eps^{-2}\cr
 \end{matrix}
 \right)                & t\leq\eps.
 \end{cases}
 $$
 We easily get $f^\pm_\eps\in\ca^\car$, $f^-_\eps\leq f\leq f^+_\eps$ 
 and $\tau(f^+_\eps-f^-_\eps)=2\eps+\eps^3$, namely $f$ is separating 
 for $\cu$-cuts in $\ca^\car$, hence $f\in\ca^\cu$.  This shows at 
 once that $\ca^\car$ is strictly smaller then $\ca^\cu$ and that the 
 latter is not an algebra, since, if $g$ is not continuous in $0$, 
 $(f^2)_{11}$ is not continuous in $0$ hence does not belong to 
 $\ca^\cu$.  The rest of the proof follows with analogous arguments.
 \end{proof}


\subsection{Densely defined trace on a non-unital C$^*$-algebra}

 The construction considered in this subsection corresponds to the most 
 general case described in the abelian setting, namely the case of a 
 Radon measure on a non-compact space. There the local structure, 
 namely the structure given by the compact support functions, plays a 
 crucial role, and the same will be in our construction. 
 
 In the classical case, a function is said Riemann integrable if it has 
 compact (or $\t$-finite) support, and is continuous but for a null set. 
 Given a C$^*$-algebra with a semicontinuous densely defined trace, we define
 a suitable analogue of this family, and show that its norm closure
 coincides with the $\car$-closure of the given C$^*$-algebra. Therefore on
 such elements, the integral may be defined in two equivalent ways. Either as
 a separating  element, as  explained before, or as a limit (when it exists)
 of the $\t$-finite restrictions (see below). 

 In the nonabelian setting, the local structure is determined by a suitable  
 family of projections which has to be closed  w.r.t. the ``$\vee$'' operation.
 We present here a family of projections with these properties; we shall show
 at the end of this section that any other reasonable family would have
 produced the same result.

 \bigskip

 In this subsection $\ca$ denotes a non unital C$^*$-algebra acting on 
 a Hilbert space $\ch$, equipped with a semicontinuous densely defined 
 trace $\t$.

 \begin{Lemma}\label{Lemma:finite}
	 Any projection in $\ca^\car$ has finite trace.
 \end{Lemma}
 \begin{proof} 
	 Let $e$ be a projection in $\ca^\car$, $a\in\ca$ such that $e\leq 
	 a$, and $a_n\in\ca_+$ with finite trace with $a_n\to a$ in norm.  
	 Then $ea_ne\to eae$, hence there exists an $n$ such that 
	 $ea_ne\geq eae-1/2e\geq1/2e$, from which the thesis follows.
 \end{proof}

 \begin{Dfn} 
	 Let us denote by $\ce$ the minimal family of projections in 
	 $\ca''$ containing all the Riemann measurable spectral 
	 projections of the selfadjoint elements in $\ca$, and closed 
	 w.r.t. the ``$\vee$'' operation.
 \end{Dfn}

 \begin{Lemma}\label{Lemma:E} 
	 The family $\ce$ is contained in $\ca^\car$.
 \end{Lemma}
 \begin{proof} 
	 Let us consider the family $\ce_1$ of projections $e\in\ar$ which 
	 are separating for $\car$-cuts $\{a^\pm_\eps\}_{\eps>0}$ such 
	 that $a^-_\eps\geq0$ and $\t(\supp (a^+_\eps)-e)<\eps$.  The 
	 $\t$-finite Riemann measurable spectral projections of a 
	 selfadjoint element $x\in\ca$ belong to $\ce_1$ by Lemma 
	 \ref{Lemma:abelian3} applied to $(\s(x),\m_x)$. \\
	 Now we show that $\ce_1$ is closed w.r.t. the $\vee$ operation.  
	 Let $e,f$ be in $\ce_1$, with $\{a^\pm_\eps\}_{\eps>0}$, 
	 $\{b^\pm_\eps\}_{\eps>0}$ the corresponding cuts.  Then 
	 $\t(e-\supp(a^-_\eps))\leq\eps$ and the same for $f$, hence 
	 $\t(e\vee f-\supp(a^-_\eps)\vee\supp(b^-_\eps))\leq2\eps$.  Since 
	 for any positive element $y$ with $\t$-finite support we have 
	 $\t(\supp(y)-y^{1/n})\to0$ and if $y_1$, $y_2$ are positive 
	 elements then $\supp(y_1)\vee\supp(y_2)=\supp (y_1+y_2)$, we get 
	 an $n$ such that $\t(e\vee 
	 f-((a^-_\eps+b^-_\eps)/2)^{1/n})\leq3\eps$.  As for the 
	 approximation from above, we observe that 
	 $\t(\supp(a^+_\eps)\vee\supp(b^+_\eps)-e\vee f)\leq2\eps$, and 
	 there exists an $n$ such that 
	 $\t((a^+_\eps+b^+_\eps)^{1/n}-\supp(a^+_\eps+b^+_\eps))\leq\eps$, 
	 hence $\t((a^+_\eps+b^+_\eps)^{1/n}-e\vee f)\leq3\eps$.  This 
	 shows that $e\vee f\in\ce_1$.  Since $\ce_1$ is contained in 
	 $\ca^\car$, the thesis follows by the minimality of $\ce$.
 \end{proof}

 Let $e$ be a projection in $\ca^\car$, and consider the C$^*$-algebra 
 $\ca_e:=\{a\in\ca:ae=ea=a\}$.  Since $\ca^\car$ is a bimodule and 
 contains $e$, we have $\ca_e+\bc e\subset\ca^\car$, hence $(\ca_e+\bc 
 e)^\car\subset\ca^\cu$.  By Lemma~\ref{Lemma:finite}, projections in 
 $\ca^\car$ are $\t$-finite.  Therefore, by the results of the 
 previous subsection, $(\ca_e+\bc e)^\car$ is an $\car$-closed 
 C$^*$-algebra in $\ca''_e$.  By Lemma~\ref{Lemma:contfunctcalc} we 
 conclude that $(\ca_e+\bc e)^\car$ is indeed contained in $\ca^\car$.  
 We may therefore consider the minimal $\car$-closed C$^*$-algebra 
 $\car_e$ verifying $\ca_e\subseteq\car_e\subset\ca^\car$.  By 
 definition, the map $e\to\car_e$ is order preserving, therefore the 
 union $\cup_{e\in\ce}\car_e$ is a $^*$-algebra, since $(\ce,\vee)$ is 
 a directed set.  We shall denote by $\car_\ce$ the C$^*$-algebra 
 given by the norm closure of $\cup_{e\in\ce}\car_e$.  Then the 
 following holds.

 \begin{Thm}\label{Thm:equal} 
	 Let $(\ca,\t)$ be a C$^*$-algebra with a semicontinuous densely 
	 defined trace, acting on a Hilbert space $\ch$.  Then 
	 $\car_\ce=\ca^\car$.
 \end{Thm}

 First we observe that such result immediately implies Theorem 
 \ref{maintheorem} in the densely defined case, by Lemma 
 \ref{keylemma}.  Before proving the Theorem, we need a simple Lemma.

 \begin{Lemma}\label{simplelemma} 
	 Let $e$ be a projection on $\ch$, and $x$ a selfadjoint operator 
	 s.t. $0\leq x\leq \a e+\b e^\perp$, $\a,\b \geq 0$.  Then 
	 $\|e^\perp xe\|\leq \sqrt{\a\b}$
 \end{Lemma}
 \begin{proof} 
	 Let $\eta_1$ be a unit vector in $e\ch$, $\eta_2$ a unit vector 
	 in $e^\perp\ch$, $\th\in(0,\pi/2)$.  Then
	 \begin{align*}
		 2\sin\th &\cos\th Re(\eta_1,x\eta_2) =\\
		 &=((\cos\th\ \eta_1+\sin\th\ \eta_2),
		 (e^\perp xe+exe^\perp)(\cos\th\ \eta_1+\sin\th\ \eta_2))\\
		 &\leq (\a\ \cos^2\th  + \b\ \sin^2\th)
	 \end{align*}
	 which gives $Re(\eta_1,x\eta_2)\leq \a\ \cot\th +\b\ \tan\th$.  
	 Taking the supremum of the left hand side w.r.t. $\eta_1$, 
	 $\eta_2$ and the infimum of the right hand side w.r.t. 
	 $\th\in(0,\pi/2)$ we immediately get the thesis.
 \end{proof}

 \begin{proof} {\it (of Theorem \ref{Thm:equal}).}
 $(\car_\ce\subset\ca^\car)$.\quad By the results of the preceding 
 subsection, it follows that $\cup_{e\in\ce}\car_e$ is a vector space 
 in $\at$, hence its norm closure $\car_\ce$ is contained in $\ar$ by 
 Lemma~\ref{Lemma:normclosure}.

 \medskip\noindent $(\ca^\car\subset\car_\ce)$.\quad Let 
 $x\in\ca^\car_\sa$, and $\{a^\pm_\eps\}_{\eps>0}$ the corresponding 
 $\car$-cut in $\ca$.  We observe that it is not restrictive to assume 
 $x\geq 0$, possibly replacing $x$ with $x-a^-_1$.  Then we set 
 $x_\d:=\phi_\d x\phi_\d$, where $\phi_\d:=\phi_\d(a^+_1)$ denotes the 
 mollified spectral projection defined in Lemma \ref{Lemma:phieps}, and 
 observe that $\{\phi_\d a^\pm_\eps\phi_\d\}_{\eps>0}$ gives an 
 $\car$-cut for $x_\d$.  Since $\phi_\d a^\pm_\eps\phi_\d\in\ca_{e_{\d}}$, 
 with $e_{\d}=\chi_{[\d/2,\infty)}(a^+_1)$, we have 
 $x_\d\in\cup_{e\in\ce}\car_e$, for 
 any $\d>0$.  The theorem is proved if we show that $x_\d\to x$ when 
 $\d\to0$.  Indeed, %
 \begin{align*}
 \|x-x_\d\| 
 &\leq\|(1-\phi_\d)x(1-\phi_\d)\|+2\|(1-\phi_\d)x\|\leq3\|e^\perp_\d 
 x\|\\
 &\leq3\|e^\perp_\d x e^\perp_\d\|+3\|e^\perp_\d x e_\d\|
 \leq3\frac{\d}{2}+6\sqrt{\frac{\d}{2}}\sqrt{\|a_{1}^{+}\|},
 \end{align*}
 where the last inequality relies on the estimates $\|e^\perp_\d x 
 e^\perp_\d\|\leq\|e^\perp_\d a^+_1 e^\perp_\d\|\leq\frac{\d}{2}$ and $x\leq 
 a^+_1\leq\frac{\d}{2} e_\d^\perp+\|a^+_1\| e_\d$, together with 
 Lemma~\ref{simplelemma}.
 \end{proof}
 
 We say that $a\in\ar$ is $\car$-summable, if at least one of $a^{+}$ 
 or $a^{-}$ is $\t$-finite.  Then we may prove
 
 \begin{Prop} 
	 $a$ is $\car$-summable iff there exists $\lim_{e\in\ar} \t(eae)$.  
	 When $a$ is positive or $\t$-finite, the above limit coincides 
	 with $\t(a)$.  
 \end{Prop}
 
 As explained before, $\t$ coincides with the pull-back via $\r_{\pi}$ 
 of the trace on the GNS representation, hence it is a semicontinuous 
 semifinite trace on $\ar$.
 
 \begin{Prop}\label{Thm:ddmain} 
	 If $\t$ is densely defined on $\ca$, then $(\ar,\t)$ is a Riemann 
	 algebra with a semicontinuous, densely defined trace.
 \end{Prop}
 \begin{proof}
	 $\t$ is densely defined on each $\car_{e}$, hence the thesis follows.
 \end{proof}

 \begin{rem} 
	 Previous Theorem shows that the algebra $\car_\ce$ may be defined 
	 independently of the family $\ce$.  We note that, replacing 
	 $\ce:=\ce_0$ with other reasonable families, the same algebra is 
	 obtained.  Indeed we may consider the set $\ce_1$ considered in 
	 the proof of Lemma~\ref{Lemma:E}, the set $\ce_2$ of all 
	 projections in $\ca^\car$, the set $\ce_3$ of all compact support 
	 projections in $\ca''$, namely projections which are majorized by 
	 an element in $\ca$ (cf.  \cite{Ped}), or the set $\ce_4$ of all 
	 $\t$-finite projections in $\ca''$.  It is easy to see that 
	 $\ce_0\subseteq\ce_1\subseteq\ce_2\subseteq\ce_3\subseteq\ce_4$, 
	 and all such sets are directed.  If $e$ is a $\t$-finite 
	 projection, $\car_e$ can be defined as the minimal Riemann 
	 algebra containing $\ca_e$, the set of such algebras being non 
	 empty, containing at least $\ca^\car$.  Then for any family 
	 $\ce_i$, we may define the C$^*$-algebra 
	 $\car_i:=\ov{\cup_{e\in\ce_i}\car_e}$, and obtain 
	 $\ca^\car=\car_0\subseteq\car_1\subseteq\car_2\subseteq\car_3\subseteq\car_4$. 
	  Since on the other hand $\car_4$ is contained in $\ca^\car$, we 
	 have proved that all these algebras coincide.  \\
	 Let $e$ be a $\t$-finite projection in $\ca''$.  We say that an 
	 $\car$-cut is $e$-supported if all elements in $A^-$ and $A^+$ 
	 belong to $\ca_e$.  Then we proved that, for any $i=0,1,2,3,4$, 
	 $\ca^\car$ is the norm closure of the Dedekind completion of 
	 $\ca$ w.r.t. $\ce_i$-supported $\car$-cuts.
 \end{rem}

 \subsection{Semifinite non-densely defined trace}\label{subsec:genunimeas}

 In the rest of the section we consider semifinite non-densely defined 
 traces, namely pairs $(\ca,\t)$ for which the closure $\ca_0$ of the 
 domain of the trace is a proper ideal in $\ca$.  
 
 \begin{Prop}\label{Thm:bimodule} 
	 The $\car$-closure $\ca^{\car}_{0}$ of $\ca_0$ is a bimodule on 
	 $\ca$.
 \end{Prop}
 \begin{proof}
	 Let $x\in (\ar_0)_{\sa}$, then there is an $\car$-cut 
	 $\{a^{\pm}_{\eps}\}$ in $\ca_{0}$ for which $x$ is a separating 
	 element.  Let $a\in\ca$, then $a^{*}xa$ is a separating element 
	 for the $\car$-cut $\{a^{*}a^{\pm}_{\eps}a\}$ in $\ca_{0}$, so 
	 that $a^{*}xa\in(\ar_0)_{\sa}$.  Now we have to show that if 
	 $x\in(\ar_0)_{\sa}$, $y\in\ar_\sa$, then $xy+yx$ and 
	 $i(xy-yx)\in(\ar_0)_{\sa}$.  Because of equations 
	 (\ref{eqn:bimo}) and (\ref{eqn:bimobis}) the thesis follows by 
	 the previous result.
 \end{proof}

 \begin{Thm}\label{Lemma:unimeaslast} 
	 The set $\ca^{\car}_{0}+\ca$ is a C$^*$-algebra, therefore 
	 coincides with $\ar$ and Theorem \ref{maintheorem} follows for 
	 $\ca$.
 \end{Thm}
 \begin{proof} 
	 Since $\ar_0$ is a $^*$-bimodule on $\ca$, $\ar_0+\ca$ is a 
	 $^*$-algebra.  Moreover, $\ar_0$ is a closed ideal in the closure 
	 of $\ar_0+\ca$, hence, by Corollary~1.5.8 in \cite{Ped} 
	 $\ar_0+\ca$ is a C$^*$-algebra.  Finally we show that $\ar\subset 
	 \ar_0+\ca$, and the thesis will follow from Lemma \ref{keylemma}.  
	 Indeed, if $x\in\ar_{\sa}$, it is a separating element for an 
	 $\car$-cut $\{a_{\eps}^{\pm}\}$ in $\ca_{\sa}$.  Then 
	 $x-a_{1}^{-}$ is a separating element for the $\car$-cut 
	 $\{a_{\eps}^{\pm}-a_{1}^{-}\}$ in $(\ca_{0})_{\sa}$, because 
	 $-\eps \leq \t(a_{\eps}^{-}-a_{1}^{-}) <1$, and $0\leq 
	 \t(a_{\eps}^{+}-a_{1}^{-})< 1+\eps$.
 \end{proof}

 \begin{Lemma}\label{Lemma:finite1} 
	 If $x\in (\ca_0^{\car}+\ca)_+$ and $\t(x)$ is finite, then $x\in 
	 \ca_0^{\car}$.
 \end{Lemma}
 \begin{proof} 
	 Let $x=r+a$, $r\in(\ar_0)_\sa$, $a\in\ca_\sa$.  Since $\t$ is 
	 densely defined on $\ar_0$ by Theorem~\ref{Thm:ddmain}, 
	 $r\in\ov{\cai(\t|_{\ar_0})}\subset\ov{\cai(\t|_{\ar_0+\ca})}$, 
	 hence $a\in\ov{\cai(\t|_{\ar_0+\ca})}$.  By 
	 Lemma~\ref{Lemma:phieps} we get $\t(\phi_\eps(|a|))<\infty$ for 
	 any $\eps>0$, hence $a\in\ov{\cai(\t|_\ca)}=\ca_0$, and the 
	 thesis follows.
 \end{proof}

 \begin{Cor}\label{Cor:structure} 
	 The $\car$-closure of a C$^*$-algebra $\ca$ with a semicontinuous 
	 semifinite trace coincides with $\ca_0^{\car}+\ca$.  The closure 
	 of the domain of $\t|_{\ar}$ coincides with $\ca_0^{\car}$.
 \end{Cor}
 \begin{proof}   
	 The first statement follows by Theorem~\ref{Lemma:unimeaslast}.  
	 The second statement is an immediate consequence of 
	 Lemma~\ref{Lemma:finite1}.
 \end{proof}

 \section{Unbounded Riemann integration}
 
 In this Section we give the construction of the unbounded Riemann 
 measurable elements.  They form a family which is closed under the 
 $\ar-\ar$ $^{*}$-bimodule operations, even if the $^{*}$-bimodule 
 properties only hold $\t- a.e.$ The product of two unbounded Riemann 
 measurable elements is not Riemann measurable in general, but, if 
 $S,\,T$ are Riemann measurable, there is $\widetilde T$, $\t$-a.e. 
 equivalent to $T$, s.t. $S\widetilde T$ is Riemann measurable.
 
 As a matter of fact, this construction will be performed on a general 
 Riemann algebra with a semicontinuous semifinite trace, and then 
 particularized to the $\car$-closure of a C$^{*}$-algebra.  We begin 
 with technical results on unbounded operators and operations on them, 
 needed in the sequel.

 Let us denote by $\cl(\ch)$ the set of linear operators on $\ch$, 
 neither necessarily bounded nor closed, with 
 $\cc(\ch)\subset\cl(\ch)$ the set of closed, densely defined 
 operators, and with $\cb(\ch)$, as usual, the C$^{*}$-algebra of 
 bounded operators.  Let $T\in\cl(\ch)$, then $T$ is densely defined 
 as an operator from $\ck:=\ov{\cd(T)}$ to $\ch$, and we can take its 
 adjoint $T^{*}: \cd(T^{*})\subset \ch \to \ck$.  Set
 $$
 T^{+}\xi :=
 \begin{cases}
 T^{*}\xi & \xi\in\cd(T^{*}) \\
 0 & \xi\in\cd(T^{*})^{\perp},
 \end{cases}
 $$
 extended by linearity.
 
 \begin{Lemma} 
	 $T^{+}\in\cc(\ch)$.
 \end{Lemma}
 \begin{proof}
	 Let $\{\xi_{n}\}\subset \cd(T^{+})$, $\xi_{n}\to\xi$, 
	 $T^{+}\xi_{n}\to\eta$, and let $p$ be the projection onto 
	 $\ov{\cd(T^{*})}$.  Then $p\xi_{n}\in\cd(T^{*})$, $p\xi_{n}\to 
	 p\xi$ and $p^{\perp}\xi_{n}\to p^{\perp}\xi\in 
	 \cd(T^{*})^{\perp}$, so that $T^{*}p\xi_{n} = T^{+}p\xi_{n}+ 
	 T^{+}p^{\perp}\xi_{n} = T^{+}\xi_{n}\to\eta$.  As $T^{*}$ is 
	 closed, we get $p\xi\in\cd(T^{*})$ and $T^{*}p\xi=\eta$, so that 
	 $\xi=p\xi+p^{\perp}\xi\in\cd(T^{+})$ and 
	 $T^{+}\xi=T^{+}p\xi=T^{*}p\xi=\eta$.
 \end{proof}
 
 \begin{Dfn} 
	 For any linear operator $T$ its {\it natural extension} is the 
	 closed operator $T^{\nat}:=(T^{+})^{*}$.  Introduce the set of 
	 {\it locally bounded} operators, w.r.t. a projection $e$, 
	 $\cl_{0}(e):= \{T\in\cl(\ch):$ there is an increasing sequence of 
	 projections $e_{n}\inc e$ s.t. $\ch_{0}:=\cup 
	 e_{n}\ch\subset\cd(T)\cap \cd(T^{*})$, $eTe_{n}$, 
	 $e_{n}Te\in\cb(\ch)\}$.
 \end{Dfn}
 
 We want to show that the natural extension of a locally bounded 
 linear operator is locally bounded as well. For the elementary 
 rules of calculus with unbounded operators we refer the reader to 
 \cite{RN} (in particular chap. 8).
  
 \begin{Prop} \label{naturale}
	 Let $T$ be a locally bounded operator w.r.t. $e$.  Then $T^{+}$ 
	 and $T^{\nat}$ are locally bounded w.r.t. $e$, too.  Moreover, 
	 setting $T_{m}:= \ov{(eTe)|_{\ch_{0}}}$ and $T_{M}:= 
	 ((eTe)^{*}|_{\ch_{0}})^{*}$, $T_{m}$, $T_{M}$ are closed 
	 operators s.t. $(T_{M})^{*}=(T^{*})_{m}$, and $T_{m} \subset 
	 \ov{eT^{\nat}e} \subset T_{M}$.
 \end{Prop}
 \begin{proof}
	 We divide it in steps.
	 \itm{i} $\cd(Te)$ is dense in $\ch$. \\
	 As $\cd(Te) =\{ \xi\in\ch : e\xi\in\cd(T)\} \equiv 
	 \cd(T)\cap e\ch \oplus e^{\perp}\ch$ and $\cd(T)\cap 
	 e\ch\supset\ch_{0}$ is dense in $e\ch$, we get the thesis.
	 
     \itm{ii} $eTe$ is closable and densely 
     defined. \\
	 Indeed, $eTe$ is densely defined because of $(i)$.  Let 
	 $\{\xi_{n}\}\subset\cd(eTe)$, $\xi_{n}\to0$, $eTe\xi_{n}\to\eta$.  
	 Then $e_{m}eTe\xi_{n}=e_{m}Te\xi_{n}\to0$, because $e_{m}Te$ is 
	 bounded.  On the other hand $e_{m}eTe\xi_{n}\to e_{m}\eta$, so 
	 that $e_{m}\eta=0$ for all $m\in\bn$.  Therefore $\eta=e\eta=0$.
	 
	 \itm{iii} $(eTe)^{*}|_{\ch_{0}}= (eT^{+}e)|_{\ch_{0}}$. \\
	 Let us set $\ck:=\ov{\cd(T)}$ and observe that $eT^{*}\equiv 
	 e^{*}T^{*}:\cd(T^{*})\subset\ch\to\ck\hookrightarrow\ch$ is 
	 extended by $(Te)^{*}$, so that 
	 $(e_{n}Te)^{*}=(Te)^{*}e_{n}\supset eT^{*}e_{n}$, where the 
	 equality holds because $Te$ is densely defined.  Therefore 
	 $eT^{+}e_{n}\subset (e_{n}Te)^{*}$, so that equality holds 
	 because $eT^{+}e_{n}$ is bounded and everywhere defined.  Finally 
	 $(eTe)^{*}e_{n} = (Te)^{*}ee_{n} = (Te)^{*}e_{n} = eT^{+}e_{n}$ 
	 which implies the thesis.
	 
	 \itm{iv} $e_{n}T^{+}e \subset (eTe_{n})^{*}$. 
	 \\
	 Let $\xi\in\cd(T^{*})\cap e\ch$, $\eta\in\ch$, then
	 \begin{align*}
		 (\eta,e_{n}T^{+}e\xi) & = (\eta,e_{n}T^{+}\xi) = (\eta,e_{n}T^{*}\xi) 
		 = (e_{n}\eta,T^{*}\xi) \\
		 & = (Te_{n}\eta,\xi) = (Te_{n}\eta,e\xi) = (eTe_{n}\eta,\xi) 
		 = (\eta,(eTe_{n})^{*}\xi)
	 \end{align*}
	 because $e_{n}\eta\in\cd(T)$ and $eTe_{n}$ is bounded.
	 
	 \itm{v} $eTe|_{\ch_{0}} = eT^{\nat}e|_{\ch_{0}}$. \\
	 As $e_{n}T^{+}e \subset (eTe_{n})^{*} \imply eTe_{n} \subset 
	 (e_{n}T^{+}e)^{*}$, we get $eTe_{n} = (e_{n}T^{+}e)^{*} \supset 
	 e(e_{n}T^{+})^{*} = eT^{\nat}e_{n}$.  But $e_{n}\ch \subset 
	 \cd(T^{\nat})$, implies $eT^{\nat}e_{n}\in\cb(\ch)$, and the 
	 thesis follows.

	 Therefore, the first statement of the proposition follows from 
	 $eT^{+}e_{n}= (e_{n}Te)^{*}$ and $e_{n}T^{+}e\subset 
	 (eTe_{n})^{*}$.  From $(ii)$ it follows that $eT^{+}e$ is 
	 closable, therefore $\ov{eT^{+}e} \supset 
	 \ov{(eT^{+}e)|_{\ch_{0}}} = \ov{(eT^{*}e)|_{\ch_{0}}} = 
	 T_{M}^{*}$, so that $T_{M} \supset (eT^{+}e)^{*} \supset 
	 e(eT^{+})^{*} = eT^{\nat}e$.  Finally $\ov{eT^{\nat}e} \supset 
	 \ov{eT^{\nat}e|_{\ch_{0}}} = \ov{eTe|_{\ch_{0}}} = T_{m}$.
 \end{proof}
  
 \bigskip
 

 Let $\car\subset\cb(\ch)$ be a Riemann algebra w.r.t.  a 
 semicontinuous semifinite trace $\t$.  Inspired by Christensen 
 \cite{Christensen} we now introduce the set of essentially 
 $\t$-measurable operators (a subset of the locally bounded ones).  
 Recall that $T\aff\car''$ stands for $u'T\subset Tu'$ for any unitary 
 operator $u'\in\car'$.
 
 \begin{Dfn} 
	 A sequence $\{e_n\}$ of projections in $\cb(\ch)$ is called a 
	 {\it Strongly Dense Domain} (SDD) w.r.t. $(\car,\t)$, if 
	 $e_n^{\perp}\in\car$ is $\tau$-finite, $\tau(e_n^{\perp})\to0$.  
	 Then $e^{\perp}:=\inf_ne^{\perp}_n \in\car$, because of 
	 $\car$-closedness.  If $T\aff\car''$, we say $\{e_n\}$ is a SDD 
	 for $T$ if $\ch_{0}:=\cup e_{n}\ch\subset\cd(T)\cap \cd(T^{*})$ 
	 and $eTe_{n}$, $e_{n}Te\in\car$.  Let us introduce the set of 
	 {\it essentially $\t$-measurable} operators $\rio := \{ 
	 T\aff\car'' :$ there is a SDD for $T \}$, and the set of {\it 
	 $\t$-measurable} operators $\ri:= \rio\cap\cc(\ch)$.  Let $S,\ 
	 T\aff\car''$, we say that $S=T$ almost everywhere if there is a 
	 common SDD $\{e_{n}\}$ for $S$ and $T$ s.t. $eSe|_{\ch_{0}} = 
	 eTe|_{\ch_{0}}$.
 \end{Dfn}
 
 \begin{rem} 
	 In order for $\ri$ to be larger than $\car$, one needs SDD's for 
	 which $\tau(e_n^{\perp})>0$ for any $n$. This is not always the 
	 case. For example  the  compact operators with the usual trace form 
	 a Riemann algebra for which $\ri=\car$.
 \end{rem}

 \begin{Lemma} \label{l:projinC} 
	 Let $\car$ be a C$^{*}$-algebra, $e,\, f$ projections in $\car$.  
	 Then $e\vee f$, $e\wedge f\in\car$.
 \end{Lemma}
 \begin{proof}
	 As $e\wedge f \leq e+f \leq 2 (e\wedge f)$ and $t\to t^{1/n}$ is 
	 operator-increasing, it follows $e\wedge f \leq (e+f)^{1/n} \leq 
	 2^{1/n} (e\wedge f)$.  Therefore $(e+f)^{1/n} \to e\wedge f$ and 
	 $e\wedge f\in\car$.  Besides, in the unitalization of $\car$, 
	 $e\vee f = 1- (1-e)\wedge(1-f) = \lim_{n} 1 - (2-e-f)^{1/n} = 
	 \lim_{n} 1 - 2^{1/n} \sum_{k=0}^{\infty} (-1)^{k} \binom{1/n}{k} 
	 \left(\frac{e+f}2\right)^{1/n} = \lim_{n} e_{n}$, where $e_{n}:= 
	 2^{1/n} \sum_{k=1}^{\infty} (-1)^{k-1} \binom{1/n}{k} 
	 \left(\frac{e+f}2\right)^{1/n}\in\car$.
 \end{proof}
 
 \begin{Lemma} \label{l:common} 
	 Let $T\in\rio$, $\{e_{n}\}$, $\{f_{n}\}$ SDD for $T$.  Then 
	 $g_{n}:= e_{n}\wedge f_{n}$ is an SDD for $T$.
 \end{Lemma}
 \begin{proof}
	 Because of Lemma \ref{l:projinC}, $g_{n}^{\perp}\in\car$, and 
	 $\t(g^{\perp}_{n})\to0$, that is $\{g_{n}\}$ is an SDD. Besides 
	 $g_{n}Tg = g_{n}e_{n}Teg = e_{n}Teg - g^{\perp}_{n}e_{n}Teg = 
	 e_{n}Te - e_{n}Teg^{\perp} - g^{\perp}_{n}e_{n}Te + 
	 g^{\perp}_{n}e_{n}Teg^{\perp} \in\car$, and the thesis follows.
 \end{proof}
 
 \begin{Prop} \label{p:naturale}
 \itm{i} $T\in\cl(\ch),\,T\aff\car'' \imply 
 T^{+},\,T^{\nat}\aff\car''$. \label{le:natisaff}	 
 \itm{ii} Equality almost everywhere is an equivalence relation.
 \itm{iii} $T\in\rio\ \imply\ T=T^{\nat}$ almost everywhere and $T^{\nat}\in\rio$.
 \end{Prop}
  \begin{proof}
     $(i)$	 Let $u'\in\car'$ be unitary, then $u'T\subset Tu'$ so that 
	 $u'T^{*}\subset T^{*}u'$, therefore $u'\cd(T^{*}) = \cd(T^{*})$.  
	 Moreover for $\eta\in\cd(T^{*})^{\perp} = u'\cd(T^{*})^{\perp}$ 
	 we get $T^{+}u'\eta = 0 = u'T^{+}\eta$.  In all 
	 $T^{+}\aff\car''$, so that $T^{\nat}\aff\car''$ as well. \\
	 $(ii)$ is obvious. \\
	 $(iii)$ As a SDD for $T$ is also a SDD for $T^{\nat}$, the thesis 
	 follows from $(i)$ and $(v)$ in the proof of Proposition \ref{naturale}.
 \end{proof}

 In particular we have just proved that the equivalence class of an 
 operator in $\rio$ contains a closed operator, hence 
 $\rio/\hskip-1.2mm\sim =\ri/\hskip-1.2mm\sim$.

 \begin{Dfn} 
	 For $S,\ T\in\ri$, $a\in\car$ define $S\oplus T:= (S+T)^{\nat}$, 
	 $a\odot T:= (aT)^{\nat}$, $T\odot a :=(Ta)^{\nat}$.
 \end{Dfn}
 
 \begin{Thm} \label{thm:unbounded-bimodule}
	 $\ri$ is an almost everywhere $^{*}$-bimodule over $\car$, w.r.t. 
	 strong sense operations, namely $S,\,T\in\ri,\,a\in\car \imply 
	 S\oplus T\in\ri,\, a\odot T,\,T\odot a\in\ri$, and the bimodule 
	 properties hold almost everywhere.
 \end{Thm}
 \begin{proof}
	 Let us first prove
	 \itm{i} $S,\,T\in\rio \imply S+T\in\rio$ and $(S+T)^{\nat} = 
	 S^{\nat}+T^{\nat}$ almost everywhere. \\
	 Assume $\{e_{n}\}$, $\{f_{n}\}$ are SDD for $S$ and $T$ 
	 respectively, and set $g_{n}:=e_{n}\wedge f_{n}$.  Then 
	 $\{g_{n}\}$ is a common SDD for $S$ and $T$, as in Lemma 
	 \ref{l:common}.  Besides $g_{n}(S+T)g = g_{n}e_{n}Seg + 
	 g_{n}f_{n}Tfg \in\car$, $g(S+T)g_{n} = geSe_{n}g_{n} + 
	 gfTf_{n}g_{n} \in\car$, as in Lemma \ref{l:common}.  Finally 
	 $g(S+T)^{\nat}g|_{\ch_{0}} = g(S+T)g|_{\ch_{0}} = gSg|_{\ch_{0}} 
	 + gTg|_{\ch_{0}} = g(S^{\nat}+T^{\nat})g|_{\ch_{0}}$. 
	 
	 Then we prove
	 \itm{ii} $T\in\rio,\,a\in\car \imply aT,\,Ta\in\rio$, and $(aT)^{\nat} = 
	 aT^{\nat}$, $(Ta)^{\nat}= T^{\nat}a$ almost everywhere. \\
	 Using $(i)$, we need only prove it for $a=u$ a 
	 unitary operator.  Set $g_{n}:=e_{n}\wedge u^{*}e_{n}u \wedge 
	 ue_{n}u^{*}$.  Then $g_{n}$ is a common SDD for $T$, $uT$, $Tu$, 
	 as in Lemma \ref{l:common}, and $g_{n}Tug = g_{n}e_{n}Tuu^{*}eug 
	 \in\car$, $gTug_{n} = geTuu^{*}e_{n}ug_{n} \in\car$, as in Lemma 
	 \ref{l:common}.  Analogously $g_{n}uTg \in\car$, $guTg_{n} 
	 \in\car$.  The last two statements are proved as in $(i)$.  
	 
	 Now the Theorem follows by  $(i)$, $(ii)$, and Proposition \ref{p:naturale}.
 \end{proof}
  
 Observe that the previous Theorem generalises results by Segal 
 \cite{Se}, indeed, if $\car$ is a von Neumann algebra and $\t$ a 
 normal semifinite faithful trace on it, equality almost everywhere 
 turns out to be equality (cf. \cite{Se} Corollary 5.1) and the two 
 notions of strong sense operations coincide.
 
  \begin{Prop}\label{Prop:quasiprod} 
	  If $S,T\in\ri$ and $\{e_n\}$ is an SDD both for $S$ and $T$ then 
	  $SeT$ and $TeS$ belong to $\ri$.  In particular, there exists 
	  $\widetilde T$ $\t$-a.e. equivalent to $T$ such that 
	  $S\widetilde T$ and $\widetilde TS$ belong to $\ri$.  
  \end{Prop}
  \begin{proof} 
	  Let $eTe_n=\sum_i\l_{ni}u_{ni}$, $e_{n}Se=\sum_i\m_{ni}v_{ni}$ 
	  be decompositions of $eTe_n$ and $e_{n}Se$ into linear 
	  combinations of unitaries in $\car$, and set 
	  $f_n:=\wedge_iu_{ni}^*e_nu_{ni}\wedge 
	  e_n\wedge_{i}v_{ni}e_nv_{ni}^{*}$.  Then $f_n\ch\subset\cd(T)$ 
	  and
	  \begin{align*}
		  eTf_n&=eTe_nf_n=\sum_i\l_{ni} u_{ni}f_n=\sum_i\l_{ni}e_nu_{ni}f_n\\
		  &=e_n(eTe_n)f_n,
	  \end{align*}
	  as a consequence $eTf_n\ch\subset\cd(S)$, then 
	  $fSeTf_n=f(eSe_n)(eTe_n)f_n\in \car$.  Repeating this argument 
	  for $f_{n}SeTf$ we show that $\{f_n\}$ is an SDD for $SeT$.  
	  Since $\widetilde T:=\ov{eTe}$ is $\t$-a.e. equivalent to $T$ 
	  the last statement follows.
  \end{proof}
 
  \bigskip
 
 From now on we denote by $\r$ the GNS representation of $\t$, and 
 with $\cam:=\r(\car)''$.  We want to extend $\r$ to a morphism of 
 $\ri$ to the $^{*}$-algebra $\mt$ of $\t$-measurable operators 
 affiliated with $\cam$ \cite{Se}.  Let us recall that the topology of 
 convergence in measure in $\mt$ is generated by the neighborhood 
 basis $\{V(\eps,\d)\}_{\eps,\d>0}$, where $V(\eps,\d) := \{ T\in\mt 
 :$ there is a projection $p\in\cam$ s.t. $\t(p^{\perp})<\d,\,\|Tp\| 
 <\eps \}$.  Let us set $T = \tlim T_{n}$ for $T_{n} \to T$ in 
 measure.
 
 \begin{Prop} \label{l:def-rho}
	 Let $T\in\ri$.
	 \itm{i} If $\{e_{n}\}$ is an SDD for $T$, then $\tlim \r(eTe_{n})$ 
	 and $\tlim \r(e_{n}Te)$ exist and are equal.
	 \itm{ii}  $\r(T) := \tlim \r(eTe_{n})$ does not depend 
	 on the SDD $\{e_{n}\}$ and belongs to $\mt$.
 \end{Prop}
 \begin{proof}
	 \itm{i}
	 The sequence $\{\r(eTe_{n})\}$ is easily seen to be Cauchy in measure.
	 Observe that from (\cite{Stn}, Theorem 3.7) $\tlim 
	 \r(e_{n}Te_{n}) = \tlim (1-\r(e^{\perp}_{n}))\r(eTe_{n}) = 
	 \tlim \r(eTe_{n})$.  Finally $\tlim \r(e_{n}Te) = \tlim \r(e_{n}Te_{n})$, 
	 because for all $\d>0$, let $k\in\bn$ be s.t. 
	 $\t(e_{k}^{\perp})<\d$, and $n>k$, then 
	 $\|[\r(e_{n}Te)-\r(e_{n}Te_{n})](1-\r(e^{\perp}_{k}))\| = 
	 \|\r(e_{n}Te_{k})-\r(e_{n}Te_{k})\|=0$.  The thesis follows. 
	 \itm{ii}
	 Let $\{f_{n}\}$ be another SDD for $T$, and set 
	 $g_{n}:=e_{n}\wedge f_{n}$.  Then $\t(g_{n}^{\perp})\to0$, and 
	 $\|(\r(eTe_{n})-\r(fTf_{n}))(1-\r(g^{\perp}_{n}))\| \leq 
	 \|\r(e-g)\r(eTe_{n})(1-\r(g^{\perp}_{n}))\| + 
	 \|\r(f-g)\r(fTf_{n})(1-\r(g^{\perp}_{n}))\| = 0$, because 
	 $\r(e^{\perp})=\r(f^{\perp})=\r(g^{\perp})=0$, as, for example, 
	 $\r(e^{\perp})=\tlim\r(e^{\perp}_{n})$.  Indeed for all $\d>0$, 
	 let $k\in\bn$ be s.t. $\t(e_{k}^{\perp})<\d$, then for $n>k$ one 
	 has $\|\r(e^{\perp}_{n})(1-\r(e^{\perp}_{k}))\|=0$.  Finally 
	 $\|[\r(eTe_{m})-\r(eTe_{n})](1-\r(e^{\perp}_{k}))\| =0$, 
	 therefore $\r(T)\in\mt$.
 \end{proof}

 Observe that $\r(\ov{eTe})=\r(T)$, because $\r(\ov{eTe}) = 
 \tlim\r(eTee_{n}) = \r(T)$.

 \begin{Thm} 
	 The map $\r:\ri\to\mt$ is a morphism of almost everywhere 
	 bimodules.  Therefore $\t\circ\r$ is a trace on $\ri$ as an 
	 almost everywhere bimodule on $\car$, extending $\t$ on $\car$.  
	 Besides, if $S,\,T\in\ri$ and $\{e_n\}$ is a common SDD, then 
	 $\r(SeT)=\r(S)\r(T)$.
 \end{Thm}
 \begin{proof} 
	 We divide it in steps.
	 \itm{i}  $T\in\ri\ \imply\ \r(T^{*})=\r(T)^{*}$. \\
	 Indeed
	 \begin{align*}
		 \r(T^{*}) & = \tlim\r(eT^{*}e_{n}) = \tlim\r(eT^{+}e_{n}) 
		 = \tlim\r(e_{n}T^{+}e) \\
		 & = \tlim\r((eTe_{n})^{*}) = \tlim\r(eTe_{n})^{*}  = \r(T)^{*}.
	 \end{align*}
	 \itm{ii} $S,\,T\in\ri\ \imply\ \r(S\oplus T)=\r(S)\oplus\r(T)$. \\
	 Indeed using Theorem 3.3 in \cite{Stn} one has
	 \begin{align*}
		 \r(S\oplus T) & = \tlim\r(g(S+T)^{\nat}g_{n}) 
		 = \tlim\r(g(S+T)g_{n}) \\
		 & = \tlim\r(gSg_{n}) \oplus \tlim\r(gTg_{n}) = \r(S) \oplus \r(T).
	 \end{align*}
	 \itm{iii}  $a\in\car$, $T\in\ri\ \imply\ \r(a\odot T)=\r(a)\odot\r(T)$, and 
	 $\r(T\odot a) = \r(T)\odot\r(a)$. \\
 	 Assume first that $a\in\car$ is unitary.  Let $f_{n}$ be an SDD for 
	 $T$, and set $g_{n}:= f_{n}\wedge af_{n}a^{*}$.  Then, by 
	 (\cite{Stn}, Theorem 3.7)
	 \begin{align*}
		 \r(a\odot T) & = \tlim\r(g(aT)^{\nat}g_{n}) 
		 = \tlim\r(gaTg_{n}) = \tlim\r(gafa^{*}aTf_{n}g_{n})\\
		 & = \tlim(1-\r(g^{\perp}))\r(a)\r(fTf_{n})(1-\r(g^{\perp}_{n})) 
		 = \r(a)\odot\r(T).
	 \end{align*}
	 The general case follows from $(i)$ and $(ii)$ in the proof of 
	 Theorem \ref{thm:unbounded-bimodule}, and $(ii)$ right above.  The proof of 
	 $\r(T\odot a) = \r(T)\odot\r(a)$ is analogous.\\
	 
	 So it remains to prove the last statement of the Theorem, which 
	 follows from
	 $$
	 \r(SeT)=\tlim\r(e_nSeTe_n)=\tlim\r(e_nSe)\tlim\r(eTe_n)=\r(S)\r(T).
	 $$
 \end{proof}
 
 \begin{rem}
	We proved in Theorem \ref{thm:unbounded-bimodule} that 
	$\ri/\hskip-1.2mm\sim$ is a $^{*}$-bimodule over $\car$, and it 
	immediately follows from the definition of $\r$ contained in 
	Proposition \ref{l:def-rho}, that if $T\in\ri$, $T=0$ $\t$-a.e., 
	then $\r(T)=0$.  Therefore we get a bimodule map from 
	$\ri/\hskip-1.2mm\sim$ to $\mt$, which is not an isomorphism, in 
	general.  More precisely, given $T\in\ri$, then $\r(T)=0$ iff for any 
	$\eps>0$, there is an SDD $\{e_{n}\}$ for $T$ s.t. $\|eTe\|<\eps$.
 \end{rem}

 \bigskip

 Finally we show how the previous construction can be applied in order 
 to extend a semicontinuous semifinite trace on a concrete 
 C$^*$-algebra to a suitable family of unbounded operators.  Let 
 $\ca\subset\cb(\ch)$ be a C$^{*}$-algebra with a semicontinuous 
 semifinite trace $\t$, and let $\ar$ be the algebra of bounded 
 Riemann measurable elements, then we call $\ov{\ar}$ the bimodule of 
 unbounded Riemann measurable elements, and the extension of $\t$ from 
 $\ar$ to $\ov{\ar}$, provided by the previous Theorem, the 
 (noncommutative) unbounded Riemann integral.
 
 \begin{Prop} \label{prop:unboundedRiemann}
	 Let $(\ca,\t)$ be a C$^{*}$-algebra with a semicontinuous 
	 semifinite trace $\t$, and $\ar$ and $\ov{\ar}$ be as above.  
	 Then $\t$ extends to a trace on $\ov{\ar}$ as an almost 
	 everywhere bimodule on $\ar$, namely $\t(uAu^{*})=\t(A)$ for any 
	 unitary operator $u\in\ar$, and any positive operator 
	 $A\in\ov{\ar}$.  Moreover, unbounded Riemann functional calculi 
	 of bounded Riemann measurable elements are Riemann measurable, 
	 namely for any $x\in\ar_{\sa}$ and 
	 $f\in\rio(\s(x)\setminus\{0\},\m_{x})$, $f(x)\in\ov{\ar}$.
 \end{Prop}
 \begin{proof} 
	 The first statement follows by the previous results in this 
	 section.  Let $f$ be in $\rio(\s(x)\setminus\{0\},\m_{x})$.  By 
	 Proposition \ref{p:4.14} in the Appendix, there exists an SDD 
	 given by the characteristic functions of Riemann measurable sets 
	 $G_{n}$ s.t. $f|_{G_{n}}$ is Riemann measurable.  As a 
	 consequence $\chi_{G_{n}}(x)$ give an SDD for $f(x)$, and the 
	 second statement follows.
 \end{proof}
 
 \begin{rem}\label{noprod} 
	 $(i)$\quad $\ov{\ar}\cap\ca''$ is not an algebra, and it is 
	 larger then $\ar$, in general.  However, given $x\in \ov{\ar}\cap 
	 \ca''$, there is a projection $p\in\ar$ s.t. $\t(p)=0$, and 
	 $p^{\perp}xp^{\perp}\in\ar$. Namely $x$ belongs to $\ar$ $\t - a.e.$ \\
	 $(ii)$\quad If $A\in\ov{\ar}$ has finite trace or is positive, 
	 its trace may be computed as $\lim_{e\in\ar}\t(eAe)$, where the 
	 projections $e$ satisfy $eAe\in\ar$.
 \end{rem}

 \section{Singular traces on C$^{*}$-algebras}

 In this Section we construct singular traces on a C$^{*}$-algebra 
 with a semicontinuous semifinite trace.  Let us first recall that, if 
 $\cam$ is a von Neumann algebra with a normal semifinite faithful 
 trace $\t$, $\mt$ the $^{*}$-algebra of $\t$-measurable operators, 
 and $T\in\mt$, its distribution function and non-increasing 
 rearrangement, the basic building blocks for the construction of 
 singular traces \cite{GI1}, are defined as follows (cf.  e.g. 
 \cite{FK,GI1})
 \begin{align*}
 \l_T(t)&:=\t(\c_{(t,+\infty)}(|T|)) \\
 \m_T(t)&:=\inf\{s\geq0: \l_T(s)\leq t\}.
 \end{align*}
 Let now $\ca$ be a C$^{*}$-algebra with a semicontinuous semifinite 
 trace $\t$ acting on a Hilbert space $\ch$.  As follows from the 
 previous Section, the GNS representation $\rho$ of $\ca$ extends to a 
 $^*$-bimodule map from the unbounded Riemann measurable operators in 
 $\ov{\ar}$ into the measurable operators of $\cam:=\rho(\ca)''$, so that we may 
 define the distribution function (and therefore the associated 
 non-increasing rearrangement) w.r.t. $\t$ of an operator 
 $T\in\ov{\ar}$ as $\l_T=\l_{\rho(T)}$, and we get $\m_T=\m_{\r(T)}$.  
 Let us observe that, if $T\in\ov{\ar}$ is a positive (unbounded) 
 continuous functional calculus of an element in $\ar$, then 
 $\c_{(t,+\infty)}(T)$ belongs to $\ar$ a.e., therefore its 
 distribution function may be defined without using the representation 
 $\r$ as $\l_T(t)=\t(\c_{(t,+\infty)}(T))$.  With these preliminaries 
 out of the way, we may carry out the construction of singular traces 
 (with respect to $\tau$) as it has been done in \cite{GI1}.

 \begin{Dfn} 
	 An operator $T\in\ov{\ar}$ is called {\it eccentric} at $0$ if 
	 either
	 $$
	 \int_0^1\mu_T(t)dt <\infty \quad {\text{and}} \quad 
	 \liminf_{t\to0}\frac{\int_0^{2t}\m_T(s)ds}{\int_0^t\m_T(s)ds} =1 
	 $$
	 or
	 $$
	 \int_0^1\mu_T(t)dt =\infty\quad{\text{and}} \quad
	 \limsup_{t\to0}\frac{\int_{2t}^1\m_T(s)ds}{\int_t^1\m_T(s)ds} =1.
	 $$
	 It is called {\it eccentric} at $\infty$ if either
	 $$
	 \int_1^\infty\mu_T(t)dt <\infty \quad {\text{and}} \quad 
	 \limsup_{t\to\infty}\frac{\int_{2t}^{\infty}\m_T(s)ds}{\int_t^\infty\m_T(s)ds} 
	 =1
	 $$
	 or
	 $$
	 \int_1^\infty\mu_T(t)dt =\infty\quad{\text{and}} \quad
	 \liminf_{t\to\infty}\frac{\int_{1}^{2t}\m_T(s)ds}{\int_1^t\m_T(s)ds} =1.
	 $$
 \end{Dfn}
 
 The following proposition trivially holds
 
 \begin{Prop} 
	 Let $(\ca,\t)$ be a C$^*$-algebra with a semicontinuous 
	 semifinite trace, $T\in\ov{\ar}$, and let $X(T)$ denote the 
	 $^*$-bimodule over $\ar$ generated by $T$ in $\ov{\ar}$, while 
	 $X(\r(T))$ denotes the $^{*}$-bimodule over $\cam:= \rho(\ca)''$ 
	 generated by $\r(T)$ in $\ov{\cam}$.  Then \item{$(i)$} $T$ is 
	 eccentric if and only if $\r(T)$ is \item{$(ii)$} 
	 $\r(X(T))\subset X(\r(T))$.
 \end{Prop}
 
 As in the case of von~Neumann algebras, with any eccentric operator 
 (at $0$ or at $\infty$) in $\ov{\ar}$ we may associate a singular 
 trace, where the word singular refers to the original trace $\t$.  
 Indeed such singular traces vanish on $\t$-finite operators, and 
 those associated to $0$-eccentric operators even vanish on all 
 bounded operators.  Of course singular traces may be described as the 
 pull-back of the singular traces on $\cam$ via the (extended) GNS 
 representation.  On the other hand, explicit formulas may be written 
 in terms of the non-increasing rearrangement.  Since Riemann 
 integration is crucial in the extension of the trace to unbounded 
 operators, we write these formulas only in case of $0$-eccentric 
 operators. Moreover this is the case occuring in Section 5.  
 
 \begin{Thm}\label{Thm:singtrac}
	 If $T\in\ov\ar$ is $0$-eccentric and $\int_0^1\mu_T(t)dt 
	 <\infty$, there exists a generalized limit $\Lim_\om$ in $0$ such 
	 that the functional
	 $$ 
	 \t_\om(A):=\Lim_\om\left(\frac{\int_0^t\m_A(s)ds}{\int_0^{t}\m_T(s)ds}\right) 
	 \quad A\in X(T)_+
	 $$ 
	 linearly extends to a singular trace on the a.e. $^*$-bimodule $X(T)$ over 
	 $\ar$ generated by $T$, where $X(T)_+$ denotes those elements whose 
	 image under $\rho$ is positive.  If $\int_0^1\mu_T(t)dt =\infty$, 
	 the previous formula should be replaced by
	 $$ 
	 \t_\om(A):=\Lim_\om\left(\frac{\int_t^1\m_A(s)ds}{\int_{t}^1\m_T(s)ds}\right), 
	 \quad A\in X(T)_+.
	 $$
	 Such traces naturally extend to traces on $X(T)+\ar$.
 \end{Thm}

 We conclude this section mentioning that the proof of Lemma~2.5 in 
 \cite{GI1} contains a gap, even though the statement is correct.  We 
 thank B. De Pagter and F. Sukochev for having noticed this gap and 
 having also furnished a correct proof.  Since the Lemma is quite 
 standard, stating that the unique positive dilation invariant 
 functional on the cone of bounded, right continuous functions with compact 
 support on $(0,+\infty)$ is the Lebesgue measure (up to a positive 
 constant), we do not include the proof here.

 \section{Novikov-Shubin invariants and singular traces}
 \label{sec:singular}
 
 In this section we show how results developed in the previous 
 sections can be applied to define and study Novikov-Shubin invariants 
 on amenable (open) manifolds with bounded geometry. More precisely we 
 assume that our manifold
 \begin{itemize}
 	 \item is a complete Riemannian manifold 
	 \item has $C^{\infty}$-bounded geometry, $i.e.$ it has positive 
	 injectivity radius, and curvature tensor bounded, with all its 
	 covariant derivatives 
	 \item is endowed with a regular exhaustion $\ck$ \cite{Roe1}, that is 
	 with an increasing sequence $\{K_{n}\}$ of compact subsets of $M$, 
	 whose union is $M$, and such that, for any $r>0$
	 $$									
	 \lim_{n\to\infty} \frac{vol(Pen^{+}(K_{n},r)) }{vol(Pen^{-}(K_{n},r))} =1,
	 $$
	 where we set $Pen^{+}(K,r):= \{x\in M: \d(x,K)\leq r\}$, and 
	 $Pen^{-}(K,r):=$ the closure of $M\setminus Pen^{+}(M\setminus K,r)$.
 \end{itemize}

 Let $M$ be as above, $F$ be a finite dimensional Hermitian vector 
 bundle over $M$, and consider the C$^{*}$-algebra $\ca(F)$ of almost local 
 operators on $L^{2}(F)$, namely the norm closure of the 
 $^{*}$-algebra of finite propagation operators, where 
 $A\in\cb(L^2(F))$ has finite propagation if there is a constant 
 $u_A>0$ s.t. for any compact subset $K$ of $M$, any $\f\in L^2(F)$, 
 $\supp\f\subset K$, we have $\supp A\f \subset Pen^+(K,u_A)$.

 \begin{Thm} \label{tutto} {\rm \cite{GI2}} 
	 \itm{i}  $\ca(F)$ contains all compact operators, 
	 \itm{ii} if $F=\La^{p} T^{*}M$ is the bundle of $p$-forms on $M$, 
	 then $f(\D_{p})\in\ca_{p}:=\ca(\La^{p} T^{*}M)$, for any $f\in 
	 C_{0}([0,\infty))$, where $\D_{p}$ is the $p$-Laplacian on $M$, 
	 \itm{iii} there exists a semicontinuous semifinite (non-densely 
	 defined) trace $Tr_{\ck}$ on $\ca(F)$, which vanishes on compact 
	 operators and, on the set of uniform operators of order $-\infty$ 
	 {\rm \cite{Roe1}}, is finite, and assumes the following form
	 $$
	 Tr_{\ck}(A) 
	 =\Lim_{\om}\frac{\int_{K_{n}}tr(a(x,x)dvol(x)}{vol(K_{n})}\ ,
	 $$
	 where $a(x,y)$ is the kernel of $A$.
 \end{Thm}

 \subsection{Novikov-Shubin numbers for open manifolds and their 
 invariance} \label{subsec:laplacian}

 In this subsection we define the Novikov-Shubin numbers for the 
 mentioned class of manifolds and prove their invariance under 
 quasi-isometries.

 Applying the results of Section 2 to $\ca_p$, we obtain the 
 C$^*$-algebra $\ca^\car_p$ with a lower-semicontinuous semifinite 
 trace, still denoted $Tr_\ck$.  Then $\chi_{[0,t)}(\D_p)$ and 
 $\chi_{[\eps,t)}(\D_p)$ belong to $\ar_p$ for almost all $t>\eps>0$, 
 by Proposition \ref{Prop:manyproj}.  Denote by $N_{p}(t) := 
 Tr_\ck(\chi_{[0,t)}(\D_p))$, $\th_{p}(t):=Tr_\ck(\e{-t\D_p})$.

 \begin{Lemma}
	$\th_p(t) = \int_{0}^{\infty}\e{-t\l} dN_p(\l)$ so that 
	$\lim_{t\to0} N_p(t)= \lim_{t\to\infty}\th_p(t)$.
 \end{Lemma}
 \begin{proof}
	If $\D=\int_{0}^{\infty}\l de(\l)$ denotes the spectral 
	decomposition, then $\e{-t\D}=\int_{0}^{\infty}\e{-t\l} de(\l)$.  
	Since the latter is defined as the norm limit of the 
	Riemann-Stieltjes sums, $\pi_p(\e{-t\D}) = 
	\int_{0}^{\infty}\e{-t\l} d\pi_p(e(\l))$, where $\pi_p$ denotes 
	the GNS representation of $\ca_p$ w.r.t the trace $Tr_\ck$.  The 
	result then follows by the normality of the trace in the GNS 
	representation.
 \end{proof}
 
 \begin{Dfn} \label{def:NS-inv} 
	 We define $b_{p}\equiv b_p(M,\ck) :=\lim_{t\to0} N_{p}(t)= 
	 \lim_{t\to\infty}\th_{p}(t)$ to be the $p$-th L$^{2}$-Betti 
	 number of the open manifold $M$ endowed with the exhaustion 
	 $\ck$.  Let us now set $N^{0}_{p}(t) := N_{p}(t) - b_{p} \equiv 
	 \lim_{\eps\to0} Tr_\ck(\chi_{[\eps,t)}(\D_p))$, and 
	 $\th^{0}_{p}(t) := \th_{p}(t) - b_{p} = \int_{0}^{\infty}\e{-t\l} 
	 dN^{0}_{p}(\l)$.  The Novikov-Shubin numbers of $(M,\ck)$ are 
	 then defined as
	 \begin{align*}
		 \a_{p}\equiv \a_{p}(M,\ck) & := 2\limsup_{t\to0} \frac{\log 
		 N^{0}_{p}(t)}{\log t}, \\
		 \underline{\a}_{p}\equiv \underline{\a}_{p}(M,\ck) & := 
		 2\liminf_{t\to0} \frac{\log N^{0}_{p}(t)}{\log t}, \\
		 \a'_{p}\equiv \a'_{p}(M,\ck) & := 2\limsup_{t\to\infty} 
		 \frac{\log \th^{0}_{p}(t)}{\log 1/t}, \\
		 \underline{\a}'_{p}\equiv \underline{\a}'_{p}(M,\ck) & := 
		 2\liminf_{t\to\infty} \frac{\log \th^{0}_{p}(t)}{\log 1/t}.
	 \end{align*}
 \end{Dfn}
 
 It follows from (\cite{GS}, Appendix) that $\underline{\a}_{p} = 
 \underline{\a}'_{p}\leq \a'_{p}\leq \a_{p}$, and 
 $\a'_{p} = \a_{p}$ if $\th_{p}^{0}(t)= O(t^{-\d})$, for 
 $t\to\infty$, or equivalently $N^{0}_{p}(t)=O(t^{\d})$, for $t\to0$.
 Observe that L$^{2}$-Betti numbers and Novikov-Shubin numbers depend 
 on the limit procedure $\om$ and the exhaustion $\ck$.
 
 \begin{rem} \label{r:remarks}
	$(a)$ The L$^{2}$-Betti numbers for amenable manifolds of bounded 
	geometry have been defined by Roe \cite{Roe1}, and it is easy to 
	show that the two definitions agree (see \cite{GI2}).  Moreover 
	Roe proved \cite{RoeBetti} that they are invariant under 
	quasi-isometries (see below).  \\
	$(b)$ If $M$ is a covering of a compact manifold $X$, L$^{2}$-Betti 
	numbers were introduced by Atiyah \cite{Atiyah} whereas 
	Novikov-Shubin numbers were introduced in \cite{NS1}.  They were 
	proved to be $\G$-homotopy invariants, where $\G:=\pi_{1}(X)$ is 
	the fundamental group of $X$, by Dodziuk \cite{Dodziuk} and 
	Gromov-Shubin \cite{GS} respectively.  \\
	$(c)$ In the case of coverings, the trace $Tr_\Gamma$ is normal on 
	the von~Neumann algebra of $\G$-invariant operators, hence 
	$\lim_{t\to0}Tr(e_{[0,t)}(\D_p))=Tr(e_{\{0\}}(\D_p))$.  In the 
	case of open manifolds there is no natural von~Neumann algebra 
	containing the bounded functional calculi of $\D_p$ on which the 
	trace $Tr_\ck$ is normal, hence the previous equality does not 
	necessarily hold.  Such phenomenon was already noticed by Roe 
	\cite{Roe2}. It has been considered by Farber 
	in \cite{Farber} in a more general context, and the 
	difference $\lim_{t\to0}Tr(e_{[0,t)}(\D_p))-Tr(e_{\{0\}}(\D_p))$ 
	has been called the {\it torsion dimension}.  We shall denote by 
	$\tordim(M,\D_{p})$ such difference.  \\
	$(d)$ Let us observe that the above definitions for $L^2$-Betti 
	numbers and Novikov-Shubin numbers coincide with the classical 
	ones in the case of amenable coverings, if one chooses the 
	exhaustion given by the F\o lner condition.  An explicit argument 
	is given in \cite{GI5}.  
 \end{rem}
 
 We prove now that Novikov-Shubin numbers are invariant under 
 quasi-isometries, where a map $\f: M\to \widetilde M$ between open 
 manifolds of \bg is a quasi-isometry \cite{RoeBetti} if $\f$ is a 
 diffeomorphism s.t.
 \begin{itemize}
	 \itm{i} there are $C_{1},\ C_{2}>0$ s.t. $C_{1}\|v\| \leq 
	 \|\f_{*}v\| \leq C_{2}\|v\|$, $v\in TM$ 
	 \itm{ii} $\nabla-\f^{*}\widetilde\nabla$ is bounded with all its 
	 covariant derivatives, where $\nabla$, $\widetilde\nabla$ are 
	 Levi-Civita connections of $M$ and $\widetilde M$.
 \end{itemize} 

 \begin{Thm} \label{thm:invariance} 
	 Let $(M,\ck)$ be an open manifold of bounded geometry with a 
	 regular exhaustion, and let $\f:M\to\widetilde M$ be a 
	 quasi-isometry.  Then $\f(\ck)$ is a regular exhaustion for 
	 $\widetilde M$, $\a_{p}(M,\ck)=\a_{p}(\widetilde M,\f(\ck))$ and 
	 the same holds for $\underline{\a}_{p}$ and $\a'_{p}$.
 \end{Thm}
 \begin{proof} 
	  Let us denote by 
	  $\Phi\in\cb(L^{2}(\La^{*}T^{*}M),L^{2}(\La^{*}T^{*}\widetilde 
	  M))$ the extension of $(\f^{-1})^{*}$.  Then 
	  $Tr_{\f(\ck)}=Tr_\ck(\Phi^{-1}\cdot\Phi)$.  Also, setting 
	  $e_{\eps,t}:=\chi_{[\eps,t)}(\D_p)$, 
	  $q_{\eta,s}:=\Phi^{-1}\chi_{[\eta,s)}(\widetilde\D_p)\Phi$, we 
	  have
	 \begin{align*}
		0 & \leq Tr_\ck(e_{\eps,t}-e_{\eps,t}q_{\eta,s}e_{\eps,t}) = 
		Tr_\ck(e_{\eps,t}(1-q_{\eta,s})e_{\eps,t}) \\
		& = Tr_\ck(e_{\eps,t}q_{0,\eta}e_{\eps,t}) + 
		Tr_\ck(e_{\eps,t}q_{s,\infty}e_{\eps,t}) \\
		& = Tr_\ck(q_{0,\eta}e_{\eps,t}q_{0,\eta}) + 
		Tr_\ck(e_{\eps,t}e_{0,t}q_{s,\infty}e_{0,t}) \\
		& \leq Tr_\ck(q_{0,\eta}e_{\eps,\infty}q_{0,\eta}) + 
		Tr_\ck(e_{\eps,t}e_{0,t}q_{s,\infty}e_{0,t}) \\
		& \leq Tr_\ck(q_{0,\eta})\|q_{0,\eta}e_{\eps,\infty}
		q_{0,\eta}\| + 
		Tr_\ck(e_{\eps,t})\|e_{0,t}q_{s,\infty}e_{0,t}\| \\
		& \leq Tr_\ck(q_{0,\eta})\ C \sqrt{\frac{\eta}{\eps}} + 
		Tr_\ck(e_{\eps,t})\ C \sqrt{\frac t s},
	\end{align*}
	where the last inequality follows from \cite{RoeBetti}.
	Then
	\begin{align*}
		Tr_\ck(q_{\eta,s}) & = Tr_\ck(e_{\eps,t}) + 
		Tr_\ck(q_{\eta,s}-e_{\eps,t}q_{\eta,s}e_{\eps,t})- 
		Tr_\ck(e_{\eps,t}-e_{\eps,t}q_{\eta,s}e_{\eps,t}) \\
		& \geq Tr_\ck(e_{\eps,t}) - Tr_\ck(q_{0,\eta})\ C 
		\sqrt{\frac{\eta}{\eps}}-Tr_\ck(e_{\eps,t})\ C\sqrt{\frac t s}.
	\end{align*}
	Now let $a>1$ and compute
	\begin{align*}
		\widetilde N^{0}(s) 
		& = \lim_{\eps\to0}Tr_\ck(q_{\eps^{a},s}) \geq 
		\lim_{\eps\to0} \left[ Tr_\ck(e_{\eps,t}) - 
		Tr_\ck(q_{0,\eps^{a}})\ C 
		\eps^{\frac{a-1}2} - Tr_\ck(e_{\eps,t})\ C \sqrt{\frac t s} 
		\right] \\
		&  = N^{0}(t)\left[ 1-C\sqrt{\frac t s} \right].
	\end{align*}
	Therefore with $\l:=4C^{2}$ we get $\widetilde N^{0}(\l t) \geq 
	\frac12 N^{0}(t)$, and exchanging the roles of $M$ and $\widetilde 
	M$, we obtain $\frac12 N^{0}(\l^{-1}t)\leq \widetilde N^{0}(t) 
	\leq 2N^{0}(\l t)$.  This means that $N^{0}$ and $\widetilde 
	N^{0}$ are dilatation-equivalent (see \cite{GS}) so that the 
	thesis follows from \cite{GS}.
 \end{proof}

 \begin{rem}
	 We have chosen Lott's normalization \cite{Lott} for the 
	 Novikov-Shubin numbers $\a_{p}(M)$, instead of the original one 
	 in \cite{NS1}.  In contrast with Lott's choice, we used the 
	 $\limsup$ in Definition \ref{def:NS-inv}.  This is motivated by 
	 our interpretation of $\a_{p}(M)$ as a dimension, as a 
	 noncommutative measure corresponds to $\a_p$ via a singular 
	 trace, according to Theorem \ref{Thm:singular}.  \\
	 In \cite{GI2} an asymptotic dimension is defined for any 
	 (noncompact) metric space.  For a suitable class of open manifolds it is 
	 shown to coincide with $\a_{0}(M)$. Therefore in this case 
	 $\a_{0}(M)$ is independent of the exhaustion and the limit procedure.
 \end{rem}

\subsection{Novikov-Shubin numbers as asymptotic spectral dimensions} 
\label{subsec:singtrac}
 In this subsection we show that Novikov-Shubin numbers can be 
 interpreted as noncommutative asymptotic dimensions.  More precisely, we prove that 
 $\a_{p}(M,\ck)$ can be expressed by a formula which is a large scale 
 analogue of Weyl's asymptotic formula for the dimension of a manifold
 $$
 \a_{p}(M,\ck) = 
 \left(\liminf_{t\to0}\frac{\log\m_{p}(t)}{\log1/t}\right)^{-1},
 $$ 
 where $\m_{p}$ refers to the operator $\D_p^{-1/2}$. When $\a_{p}$ 
 is finite non-zero, $\D_p^{-\a_{p}/2}$ gives rise to a singular 
 trace, namely there exists a type II$_{1}$ singular trace which is 
 finite nonzero on $\D_p^{-\a_{p}/2}$.

 This result, which extends the analogous result for coverings in 
 \cite{GI5}, makes essential use of the unbounded Riemann integration 
 and the theory of singular traces for C$^*$-algebras developed in 
 Sections 3 and 4.  However, since the trace we use is not normal with 
 respect to the given representation of $\ca_{p}$ on the space of 
 $L^2$-differential forms, some assumptions, like the vanishing of the 
 torsion dimension introduced in Remark \ref{r:remarks} $(c)$, are 
 needed.

 \medskip

 In the following, when the Laplacian $\D_p$ has a non trivial kernel, 
 we denote by $\D_p^{-\a}$, $\a>0$, the (unbounded) functional 
 calculus of $\D_p$ w.r.t. the function $\f_\a$ given by $\f_\a(0)=0$ 
 and $\f_\a(t)=t^{-\a}$ when $t>0$.

 \begin{Prop}\label{Lemma:tordimass} Let $M$ be an amenable open 
	 manifold.  If 
	 \itm{a} the projection $E_{p}$ onto the kernel of $\D_{p}$  is Riemann measurable, 
	 and the torsion dimension vanishes, namely $Tr_\ck(E_{p})$ is 
	 equal to $b_p$, \\
	 then $\D_{p}^{-\a}\in\ov{\ar_{p}}$ for any $\a>0$. \\
	 The vanishing of the Betti number $b_p$ implies
	 $(a)$. It is equivalent to $(a)$ if $\ker(\D_{p})$ is
	 finite-dimensional.
 \end{Prop}
 \begin{proof}
	 By hypothesis $E_{p} = \chi_{\{1\}}(\e{-\D_{p}})\in\ar_{p}$, 
	 hence $T_{p} := \chi_{[0,1)}(\e{-\D_{p}})\e{-\D_{p}} = 
	 \e{-\D_{p}} - E_{p}\in \ar_{p}$.  Then the spectral measure 
	 $\n_{p}$ associated with $T_{p}$ as in equation (\ref{eqn:muex}), 
	 is a finite measure on $[0,1]$ (see Theorem \ref{tutto}) and 
	 $\n_{p}(\{1\}) = 0$.  Therefore, the function
	 $$
	 f_{\a}(t) :=
	 \begin{cases}
	 	(-\log t)^{-\a} & t\in (0,1) \\
		0 & t=0,1.
	 \end{cases}
	 $$
	 belongs to $\rio((0,1],\n_{p})$ (see Proposition \ref{p:4.14}), 
	 and $\D_{p}^{-\a} = f_{\a}(T_{p}) = f_{\a}(\e{-\D_{p}}) 
	 \in\ov{\ar_{p}}$ by Proposition \ref{prop:unboundedRiemann}.  
	 This proves the first statement.\\
	 If $\lim_{t\to\infty} Tr_\ck(\e{-t\D_{p}})=0$, from $0\leq
	 \chi_{\{1\}}(\e{-\D_{p}}) \leq \e{-t\D_{p}}$, we have that
	 $\chi_{\{1\}}(\e{-\D_{p}})$ is a separating element for an
	 $\car$-cut in $\ca$. The last statement follows from
	 the vanishing of $Tr_{\ck}$ on compact operators (see Theorem 
	 \ref{tutto}).
 \end{proof}

 Conditions implying the vanishing of L$^{2}$-Betti numbers are given in 
 \cite{Lott2}. 
 
 \medskip

 If hypothesis $(a)$ of the previous Lemma is satisfied, we may define 
 the distribution function $\l_{p}$ and the eigenvalue function 
 $\m_{p}$ for the operator $\D_p^{-1/2}$, hence the local spectral 
 dimension as the inverse of 
 $\lim_{t\to\infty}\frac{\log\m_{p}(t)}{\log1/t}$, which may be shown 
 to coincide with the dimension of the manifold for any $p$.
 Moreover 
 
 \begin{Dfn}
	The asymptotic spectral dimension of the triple $(M,\ck,\D_p)$ is
 $$
 \left(\liminf_{t\to0}\frac{\log\m_{p}(t)}{\log1/t}\right)^{-1}.
 $$
 \end{Dfn}

 \begin{rem}
	 $(i)$ The extension of the GNS representation $\pi$ to $\ar_{p}$ 
	 does not necessarily commute with the Borel functional calculus.  
	 In particular $\chi_{\{1\}}(\pi(\e{-\D_{p}}))$ is not necessarily 
	 equal to $\pi(\chi_{\{1\}}(\e{-\D_{p}}))$.  \\
	 $(ii)$ Condition $(a)$ of the previous Proposition is equivalent to \\
	 $(a')$ $\chi_{\{1\}}(\pi(\e{-\D_{p}}))$ is Riemann integrable in 
	 the GNS representation $\pi$.  \\
	 The proof goes as follows.  \\
	 $(a)\imply (a')$.  Since the projection 
	 $E_{p}\equiv\chi_{\{0\}}(\D_p)$ is Riemann integrable and less 
	 than $\e{-t\D_p}$ for any $t$, its image in the GNS 
	 representation is Riemann integrable and less than 
	 $\pi(\e{-t\D_p})$ for any $t$.  This implies that 
	 $\pi(E_p)\leq\chi_{\{1\}}(\pi(\e{-\D_{p}}))\leq \pi(\e{-t\D_p})$ 
	 is an $\car$-cut, hence the thesis.  \\
	 $(a') \imply (a)$.  By normality of the trace in the GNS 
	 representation, $Tr_\ck(\e{-t\D_p})$ converges to 
	 $Tr_\ck(\chi_{\{1\}}(\pi(\e{-\D_{p}})))$ hence, by hypothesis, 
	 for any $\eps>0$ we may find $a_\eps\in\ca$ and $t_\eps>0$ s.t. 
	 $a_\eps\leq\chi_{\{1\}}(\pi(\e{-\D_{p}}))$ and 
	 $Tr_\ck(\e{-t_\eps\D_{p}} - a_\eps)<\eps$.  This implies 
	 $a_\eps\leq\e{-t\D_p}$ for any $t$, hence 
	 $a_\eps\leq\chi_{\{0\}}(\D_p)$, which means that 
	 $(\{a_\eps\},\{\e{-t_\eps\D_p}\})$ is an $\car$-cut for 
	 $\chi_{\{0\}}(\D_p)$, namely this projection is Riemann 
	 integrable and 
	 $Tr_\ck(\e{-t_\eps\D_p}-\chi_{\{0\}}(\D_p))\leq\eps$, i.e. the 
	 thesis.
 \end{rem}

 \begin{Thm}\label{Thm:singular} 
	 Let $(M,\ck)$ be an open manifold equipped with a regular 
	 exhaustion such that the projection on the kernel of $\D_p$ is 
	 Riemann integrable and  $\tordim(M,\D_{p})=0$.  
	 Then 
	 \itm{i} the asymptotic spectral dimension of $(M,\ck,\D_p)$ 
	 coincides with the Novikov-Shubin number $\a_p(M,\ck)$,
	 \itm{ii} if $\a_p$ is finite 
	 nonzero, then $\D_p^{-\a_p/2}$ is $0$-eccentric, therefore gives rise 
	 to a non trivial singular trace on the unbounded Riemann measurable 
	 opeators of $\ca_{p}$.
 \end{Thm}
 \begin{proof} 
	 \itm{i} By hypothesis, $e_{(0,t)}(\D_p)$ is Riemann integrable 
	 $Tr_\ck$-a.e., hence $N^0_p(t) = Tr_\ck(e_{(0,t)}(\D_p)) = 
	 Tr_\ck(e_{(t^{-1},\infty)}(\D_p^{-1})) = 
	 Tr_\ck(e_{(t^{-1/2},\infty)}(\D_p^{-1/2})) = \l_p(t^{-1/2})$.  Then
	 \begin{align}
		 \a_p &=2\limsup_{s\to0} \frac{\log N^0_{p}(s)}{\log{s}} 
		 &=2\limsup_{s\to0}\frac{\log\l_p(s^{-1/2})}{\log{s}} 
		 &=\limsup_{t\to\infty}\frac{\log\l_p(t)}{\log\frac1t}.
	 \end{align}
	 The statement follows from
	 $$
	 \liminf_{t\to0}\frac {\log\m(t)}{\log\frac1{t}} 
	 =\left(\limsup_{s\to\infty}\frac 
	 {\log\l(s)}{\log\frac1{s}}\right)^{-1}
	 $$
	 which is proved in \cite{GI5}.  
	 \itm{ii} When $0<\a_{p}<\infty$, 
	 $$
	 \liminf_{t\to0}\frac{\log \m_{\D_{p}^{-\a_{p}/2}}(t)}{\log(1/t)} = 
	 1,
	 $$
	 and this implies the eccentricity condition, as shown in 
	 \cite{GI5}.  Hence the thesis follows by 
	 Theorem~\ref{Thm:singtrac}.
 \end{proof}

 \section{Appendix}

 Here we present some more or less known results on the Riemann 
 measurable functions on a locally compact Hausdorff space $X$ which 
 are needed in the previous sections.  
 
 Let us introduce $\li := \{ f : X \to \bc: f$ is bounded and 
 $\lim_{x\to \infty}f(x)=0 \}$, and, for $\m$ an outer regular, 
 complete, positive Borel measure, $\ro := \{ f\in\li : f$ is 
 continuous but for a set of zero $\m$-measure $\}$, the set of 
 Riemann $\m$-measurable functions.  Let us observe that any 
 semicontinuous semifinite trace on $\pic{X}$ gives rise to such a 
 measure $\m$.
 
 A different description of $\ro$ is contained in the following 
 Proposition, whose proof we leave to the reader.
 
 \begin{Prop}\label{Thm:1} Let $f:X \to \br$, then the following are 
 equivalent 
 \itm{i} $f\in\ro$ 
 \itm{ii} for any $\eps>0$ there are 
 $f_{\eps}^{\pm}\in\pic{X}$ s.t.  $f_{\eps}^{-}\leq f\leq 
 f_{\eps}^{+}$ and $\int(f_{\eps}^{+}-f_{\eps}^{-})d\m<\eps$ 
 \itm{iii} there are $h\in\pic{X}$, an open subset $V$ of finite 
 $\m$-measure, and, for any $\eps>0$, $h_{\eps}^{\pm}\in\pic{V}$ 
 s.t. $h_{\eps}^{-}\leq f-h \leq h_{\eps}^{+}$ and 
 $\int(h_{\eps}^{+}-h_{\eps}^{-})d\m<\eps$.
 \end{Prop}

 We now prove some Lemmas used in Section 2.

 A measurable set $\O\subset X$ is said Riemann $\m$-measurable if its 
 characteristic function is Riemann $\m$-measurable, which is 
 equivalent to saying $\m(\partial\O)=0$.
 
 \begin{Sublemma}\label{sublemma:abelian2} Let $f$ be a positive Riemann 
 $\m$-measurable function on $X$ such that $\int f<\infty$ and set 
 $O_y=\{x\in X:f(x)> y\}$, $C_y=\{x\in X:f(x)\geq y\}$, $y\geq0$.  
 Then, for any $0<\a<\b$ there are uncountably many $y\in(\a,\b)$ for 
 which $O_{y}$ and $C_{y}$ are Riemann $\m$-measurable.
 \end{Sublemma}

 We omit the proof since it follows by standard arguments. 
 
 \begin{Lemma}\label{Lemma:abelian2} Let $f$ be a positive Riemann 
 measurable function such that $\int f d\m$ $<$ $\infty$. Then, for any 
 $\d>0$, we may find a sequence of Riemann measurable characteristic 
 functions $\chi_n$ and a sequence of positive numbers $\a_n$ such that 
 $\sum_n\a_n=\|f\|$ and $f=\sum_n\a_n\chi_n$.
 \end{Lemma}

 \begin{proof} First we construct by induction the sequences 
 $\chi_n:=\chi_{\O_n}$ and $\a_n>0$, $n\geq1$, such that, $\forall n\geq0$,
 \begin{itemize}
	 \item $f\geq\sum_1^n\a_k\chi_k$, 
	 \item $\b_{n}:=\|f-\sum_1^n\a_k\chi_k\|>0$, 
	 \item $\b_n/4\leq\a_{n+1}\leq\b_n/2$,
	 \item $\O_{n+1}:=\left\{x\in X:f(x)-\sum_1^n\a_k\chi_k(x)\geq\a_{n+1}
       \right\}$ is Riemann measurable.  
 \end{itemize}
 Indeed, given $\chi_k$, $\a_k$ for $k\leq n$, as prescribed, the 
 existence of $\a_{n+1}$ satisfying the last two properties follows by 
 Sublemma~\ref{sublemma:abelian2}, while the first two inequalities 
 (hence $\a_n>0$) follow by the definition of $\a_{n+1}$ and 
 $\O_{n+1}$.

 We now observe that, since 
 $\sup_X (f-\sum_1^n\a_k\chi_k)=\sup_{\O_{n+1}} (f-\sum_1^n\a_k\chi_k)$,
 we have 
 \begin{align*}
 \sup_{\O_{n+1}}\left(f-\sum_{k=1}^{n+1}\a_k\chi_k\right)&=  
 \sup_{\O_{n+1}}\left(f-\sum_{k=1}^n\a_k\chi_k-\a_{n+1}\right)=
 \b_n-\a_{n+1}\geq\a_{n+1}\\
 \sup_{\O_{n+1}^c}\left(f-\sum_{k=1}^{n+1}\a_k\chi_k\right)&=
 \sup_{\O_{n+1}^c}\left(f-\sum_{k=1}^n\a_k\chi_k\right)\leq\a_{n+1}.
 \end{align*}
 hence
 $$
 \b_{n+1}=
	\max\left(\sup_{\O_{n+1}}(f-\sum_{k=1}^{n+1}\a_k\chi_k), 
 	\sup_{\O^c_{n+1}}(f-\sum_{k=1}^{n+1}\a_k\chi_k)\right)
	=\b_n-\a_{n+1}\leq\frac34\b_n.
 $$
 This shows at once that $\b_n\leq(3/4)^n\|f\|\to0$, namely
 $\sum\a_n\chi_n$ converges to $f$ uniformly and
 $\sum_0^\infty\a_{n+1}=\sum_0^\infty(\b_n-\b_{n+1})=\|f\|$,
 which concludes the proof.
 \end{proof}

 \begin{Lemma}\label{Lemma:abelian3} Let  $\O\Subset X$ with 
 $\m(\partial\O)=0$ and $\m(\O)<\infty$.  Then for any $\eps>0$ $\exists f^\pm_\eps\in 
 \pic{X}$ such that $0\leq f^-_\eps\leq\chi_\O\leq f^+_\eps\leq1$, 
 $\int(f^+_\eps- f^-_\eps)d\m\leq\eps$ and 
 $\m(\supp(f^+_\eps)\setminus\O)\leq\eps$.
 \end{Lemma}

 \begin{proof} Since $\chi_\O$ is Riemann measurable, by   
 Proposition~\ref{Thm:1}, we may find $f^\pm_\eps$ satisfying all 
 the properties above, except possibly the last one.  Then, choosing a 
 continuous increasing function $\psi_\d$ on $[0,1]$ s.t.  
 $\psi_\d(t)=0$ when $t\in[0,1-\d]$ and $\psi_\d(1)=1$, we may replace 
 $f^+_\eps$ with $\psi_\d\circ f^+_\eps$.  We have
 \begin{align*}
 \m(\supp(\psi_\d\circ f^+_\eps)) &\leq\m(\{x\in X:f^+_\eps(t)>1-\d\}) \\
 &\leq\frac{1}{1-\d}\int f^+_\eps d\m\leq\frac{1}{1-\d}(\eps+\m(\O))
 \end{align*}
 from which the thesis follows.
 \end{proof}
 
 We conclude this Appendix giving a characterization, in the 
 commutative case, of the unbounded Riemann $\m$-measurable functions.
 
 \begin{Prop} \label{p:4.14} Setting $\ov\car_{0}(X,\m) := \{ f:X \to 
 \bc : f$ is $\m$-a.e.  defined and continuous, and there is a compact, Riemann 
 $\m$-measurable subset $K$ of finite $\m$-measure s.t.  
 $f|_{K^{c}}\in\li \}$, we have $\ov\car_{0}(X,\m) = \ov{(\pic{X})^\car}$, in 
 the universal atomic representation.
 \end{Prop}
 \begin{proof}
 Let $f\in\ov{(\pic{X})^\car}$.  Then there is an increasing sequence of 
 Riemann $\m$-measurable subsets $G_{n}$ s.t. $G_{n}^{c}\Subset X$, 
 $\m(G_{n}^{c})\dec0$ and $f\chi_{G_{n}}\in\ro$. Therefore 
 $\lim_{x\to\infty}f(x)=0$, and, for any $n\in\bn$ there is 
 $E_{n}\subset G_{n}$ s.t. $\m(E_{n})=0$ and $f|_{E_{n}^{c}\cap 
 G_{n}}$ is continuous. Setting $G:=\cup_{n}G_{n}$, 
 $E:=\cup_{n}E_{n}$, we get $\m(E)=\m(G^{c})=0$ and $f|_{E^{c}\cap G}$ 
 is (defined and) continuous. Finally choose an $n\in\bn$ and set 
 $K:=\ov{G_{n}^{c}}$. Then $f\in\ov\car_{0}(X,\m)$.  \\
 Let now $f\in \ov\car_{0}(X,\m)$, and let $E\subset X$ be s.t. 
 $f|_{E^{c}}$ is defined and continuous, and $\m(E)=0$. Then, because of outer 
 regularity of $\m$, there are open sets $\O_{\eps}\dec E$, with 
 $\m(\O_{\eps})<\eps$. Let us fix $\eps>0$ and set, for any $\l>0$, 
 $V_{\eps,\l}:= \{x\in K\cap \O^{c}_{\eps}: |f(x)|\leq \l \}$. Then, 
 with $M:= \sup_{K^{c}} |f|$, from Sublemma \ref{sublemma:abelian2} we 
 conclude the existence of uncountably may $\l>M$ s.t. $V_{\eps,\l}$ 
 is Riemann $\m$-measurable. For any $\eps>0$ choose 
 one such $\l\equiv\l_{\eps}$ so large that, with 
 $A_{\eps}:=V_{\eps,\l_{\eps}}^{o}$, we get $\m(K\cap A^{c}_{\eps})= 
 \m(K)-\m(V_{\eps,\l_{\eps}})<2\eps$.  Finally 
 set $G_{n}:= \cup_{k=1}^{n} A_{1/k}\cup K^{c}$, which is an 
 increasing family of Riemann $\m$-measurable subsets, s.t.  
 $\m(G_{n}^{c})\leq \m(K\cap A_{1/n}^{c})\leq \frac2{n}$, and $\sup_{G_{n}} 
 |f| = \max\{M, \max_{k\leq n} \sup_{A_{1/k}} |f|\} \leq \max_{k\leq 
 n} \l_{1/k}$.  Therefore $\{G_{n}\}$ is an SDD for $f$, and 
 $f\in\ov{(\pic{X})^\car}$.
 \end{proof}

 \begin{ack} We would like to thank Dan Burghelea, Gert Pedersen, 
 Paolo Piazza, and Laszlo Zsido for comments and suggestions.
 \end{ack} 
 


\begin{thebibliography}{99}
 
 \bibitem{Atiyah} M. F. Atiyah
{\it Elliptic operators, discrete groups and von Neumann algebras}.
Soc. Math. de France, Ast\'erisque {\bf 32--33} (1976),
43--72.

 \bibitem{ChavelFeldman} I. Chavel, E. A. Feldman. {\it Modified 
 isoperimetric constants, and large time heat diffusion in Riemannian 
 manifolds}. Duke J. Math., {\bf 64} (1991), 473--499.

 \bibitem{Christensen} E. Christensen.  {\it Non commutative 
 integration for monotone sequentially closed C$^{*}$-algebras}.  
 Math.  Scand., {\bf 31} (1972), 171--190.

 \bibitem{Co} A. Connes. {\it Non Commutative Geometry}. Academic
Press, 1994.

 \bibitem{DixmierC} J. Dixmier.  {\it C$^*$-algebras}.  North-Holland 
 Publ., Amsterdam, 1977.

 \bibitem{Dodziuk} J. Dodziuk.  {\it De Rham-Hodge theory for 
 $L^2$-cohomology of infinite coverings}.  Topology, {\bf 16} (1977), 
 157--165.

 \bibitem{FK} T. Fack, H. Kosaki.  {\it Generalized s-numbers of 
 $\t$-measurable operators}.  Pacific J. Math., {\bf 123} (1986), 269.

 \bibitem{Farber} M. Farber.  {\it Geometry of growth: approximation 
 theorems for $L^2$ invariants}.  Math.  Ann., {\bf 311} (1998), 
 335--375.
 
 \bibitem{GS} M. Gromov, M. Shubin. {\it Von Neumann spectra near
zero}. Geometric and Functional Analysis, {\bf 1} (1991),
375--404.

 \bibitem{GI1} D. Guido, T. Isola.  {\it Singular traces for 
 semifinite von~Neumann algebras}.  Journal of Functional Analysis, 
 {\bf 134} (1995), 451--485.

 \bibitem{GI5} D. Guido, T. Isola.  {\it Singular traces, dimensions, 
 and Novikov-Shubin invariants}, preprint math.OA/9907005, to appear 
 in the Proceedings of the 17th International Conference on Operator 
 Theory, Timisoara, Romania, June 1998.
 
 \bibitem{GI2} D. Guido, T. Isola.  {\it Novikov-Shubin invariants and 
 asymptotic dimensions for open manifolds}, preprint math.DG/9809040.

 \bibitem{Kadi1} R. V. Kadison.  {\it Unitary invariants for 
 representations of operator algebras}.  Ann of Math.  {\bf 66} 
 (1957), 304-379.

 \bibitem{Kaplan} S. Kaplan.  {\it The bidual of $C(X)$, I}.  
 North-Holland Publ., Amsterdam, 1985.

 \bibitem{Lott} J. Lott. {\it Heat kernels on covering spaces and
topological invariants}. J. Diff. Geom., {\bf 35} (1992),
471--510.

 \bibitem{Lott2} J. Lott. {\it On the spectrum of a finite-volume 
 negatively-curved manifold}, preprint math.DG/9908136.

 \bibitem{NS1} S. P. Novikov, M. A. Shubin. {\it Morse theory and von
Neumann {\rm II}${}_1$ factors}. Doklady Akad. Nauk SSSR, {\bf
289} (1986), 289--292.

 \bibitem{Ped} G. K. Pedersen.  {\it C$^{*}$-algebras and their 
 automorphism groups}.  Academic Press, London, 1979.

 \bibitem{RN} F. Riesz, B. Sz.  Nagy.  {\it Functional analysis}.  F. 
 Ungar Publ.  Co., New York, 1978.

 \bibitem{Roe1} J. Roe.  {\it An index theorem on open manifolds.  I}.  
 J. Diff.  Geom., {\bf 27} (1988), 87--113.

 \bibitem{Roe2} J. Roe.  {\it An index theorem on open manifolds.  
 II}.  J. Diff.  Geom., {\bf 27} (1988), 115--136.

 \bibitem{RoeBetti} J. Roe.  {\it On the quasi-isometry invariance of 
 L$^{2}$ Betti numbers}.  Duke Math.  J., {\bf 59} (1989), 765--783.

 \bibitem{Schw1} L. Schwartz.  {\it Radon measures on arbitrary 
 topological spaces and cylindrical measures}.  Oxford University 
 Press, 1973.

 \bibitem{Se} I.E. Segal.  {\it A non-commutative extension of 
 abstract integration}.  Ann.  Math.  {\bf57} (1953) 401.

 \bibitem{Stn} W.F. Stinespring.  {\it Integration theorems for gages 
 and duality for unimodular groups}.  Trans.  Am.  Math.  Soc., {\bf 
 90} (1959), 15.

 \bibitem{Tala1} M. Talagrand.  {\it Closed convex hull of set of 
 measurable functions, Riemann-measurable functions and measurability 
 of translations}.  Ann.  Inst.  Fourier, {\bf 32} (1982), 39-69.
 
 \bibitem{Takesaki} M. Takesaki. {\it Theory of operator algebras, 
 I}, Springer Verlag, Berlin, Heidelberg, New York, 1979.

 \bibitem{Varopoulos1} N. T. Varopoulos.  {\it Random walks and 
 Brownian motion on manifolds}.  Symposia Mathematica, {\bf XXIX} 
 (1987), 97--109.
 
 \bibitem{Varopoulos2} N. T. Varopoulos.  {\it Brownian motion and 
 random walks on manifolds}.  Ann.  Inst.  Fourier, Grenoble, {\bf34}, 
 2, (1984), 243--269.


\end{thebibliography}
\end{document}